\documentclass[reqno,12pt]{amsart}
\usepackage{xr}
\usepackage{graphicx}
\usepackage[usenames,dvipsnames,svgnames,table]{xcolor}
\usepackage{epstopdf}
\usepackage{booktabs}
\usepackage{amssymb,amsmath}
\usepackage{amsthm}
\usepackage{enumerate}
\newtheorem{theorem}{Theorem}[section]
\newtheorem{proposition}[theorem]{Proposition}
\newtheorem{lemma}[theorem]{Lemma}
\newtheorem{definition}[theorem]{Definition}

\newtheorem{remark}[theorem]{Remark}
\newtheorem{example}[theorem]{Example}
\numberwithin{equation}{section}

\newcommand{\R}{{\mathbb R}}
\newcommand{\CA}{{\mathcal{A}}}

\newcommand{\CV}{{\mathcal{V}}}

\newcommand{\CE}{{\mathcal{E}}}
\newcommand{\CF}{{\mathcal{F}}}

\newcommand{\la}{{\langle}}

\newcommand{\ra}{{\rangle}}

\newcommand{\q}{\mathfrak{q}}
\newcommand{\g}{\mathfrak{g}}

\newcommand{\C}{{\mathbb C}}
\newcommand{\Z}{{\mathbb Z}}
\renewcommand{\ll}{{\langle}}
\newcommand{\rr}{{\rangle}}

\newcommand{\Ch}{{\mathrm{Ch}}}
\newcommand{\ChBV}{{\mathrm{ChBV}}}

\newcommand{\Cone}{{\mathrm{Cone}}}

\newcommand{\Sym}{{\mathrm{Sym}}}

  \newcommand{\End}{{\mathrm{End}}}
\newcommand{\Tr}{{\mathrm{Tr}}}

\newcommand{\Infdex}{{\mathrm{Infdex}}}

\newcommand{\Index}{{\mathrm{Index}}}
 \newcommand{\supp}{{\mathrm{supp}}}
 \newcommand{\mult}{{\mathrm{mult}}}

\newcommand{\A}{{\mathcal{A}}}

\newcommand{\CH}{{\mathcal{H}}}

\newcommand{\CL}{{\mathcal{L}}}
\newcommand{\CN}{{\mathcal{N}}}

\newcommand{\CS}{{\mathcal{S}}}

 \renewcommand{\c}{{\mathfrak{c}}}
\newcommand{\s}{{\bf{s}}}

\renewcommand{\t}{{\mathfrak t}}
\begin{document} \title[ formula for the index]{
  Formal equivariant $\hat A$  class, splines and  multiplicities of the index of transversally elliptic operators}

\author{ Mich{\`e}le Vergne} \address{Universit\'e Denis-Diderot-Paris 7, Institut de Math\'ematiques de Jussieu,
 C.P.~7012\\ 2~place Jussieu,   F-75251 Paris~cedex~05}
\email{michele.vergne@imj-prg.fr}

\maketitle

\begin{abstract}
 Let $G$ be a  connected compact Lie group  acting on a manifold $M$ and let
 $D$ be a transversally elliptic operator on $M$.
 The multiplicity of the index of $D$ is a function  on the set $\hat G$ of irreducible representations of $G$.
Let $T$ be a maximal torus of $G$ with Lie algebra $\t$.
  We  construct a finite number of piecewise polynomial functions on $\t^*$, and
give a formula for the multiplicity in term of these functions. The main new concept is the formal equivariant $\hat A$ class. \end{abstract}

\tableofcontents

\section{Introduction} Let $G$ be a compact  Lie group acting on a manifold $M$ of dimension $d$.  As described in the monograph \cite{At}, Atiyah-Singer have associated to any
$G$-transversally elliptic symbol $\sigma$ on $M$ a   virtual trace class representation of   $G$. Let
$\Index_G(\sigma)(g)$ its trace:
$$\Index_G(\sigma)(g)=\sum_{\lambda\in \hat G} \mult_G(\sigma)(\lambda)\chi_\lambda(g).$$
Thus $\Index_G(\sigma)(g)$ is a $G$-invariant (generalized) function on $G$, and the right hand side of the above formula
is its Fourier expansion in terms of
the traces $\chi_\lambda$ of the unitary irreducible representations $V_\lambda$ of $G$.
If $D$ is a transversally elliptic operator with  principal symbol $\sigma$, we write indifferently  $\Index_G(D)$ or $\Index_G(\sigma)$, $\mult_G(D)$ or $\mult_G(\sigma)$.
The computation of $\mult_G(\sigma)(\lambda)$ is important.
For example, if $D$ is a transversally elliptic operator with principal symbol $\sigma$, the multiplicity of the trivial representation in $\Index_G(\sigma)$ is the (virtual) dimension of the space of $G$-invariant (virtual) solutions of $D$.

In this article, we will restrict ourselves to the case where $G$ is connected.

Let $\g$ be the Lie algebra  of $G$, and $\g^*$ its dual vector space.
Our aim is to construct a canonical $G$-invariant function $m_G(\sigma)$ on $\g^*$ which extends the
multiplicity function $\mult_G(\sigma)$ on $\hat G\subset \g^*/G$.

The first instance of such a relation between the multiplicity function on $\hat G$  and  functions on $\g^*$
 is the formula for the Kostant partition function in terms of derivatives of spline functions (that is piecewise polynomial functions)
\cite{DM1}, \cite{brionvergne}.
Similarly,   Heckman's result \cite{Heckman82} on branching rules relates asymptotically multiplicities  to spline functions.
   For example, if $T$ is the maximal torus of $G$,
 the asymptotic $T$ multiplicity function of the irreducible representation
$V_{k\lambda}$ of $G$  suitably normalized converges when $k$ tends to $\infty$
 to the Duistermaat-Heckman measure, a piecewise polynomial measure on $\t^*$ supported on the projection  of  $G\lambda$ on $\t^*$.
Since then, much more precise results have been obtained, in the spirit of the $[Q,R]=0$ theorem, for special elliptic symbols
\cite{Meinrenken-Sjamaar}, \cite{pep-vergne:spinc}.

Here we consider a {\bf general transversally elliptic symbol}.
 The value of our  function  $m_G(\sigma)$ at $\lambda$ is given by a double integral formula reminiscent of Witten non abelian localization formula \cite{Wittennonabelian}.
 Our present result does not provide another proof of the $[Q,R]=0$ theorems as, for general $\sigma$, we do not have a geometric interpretation of  $m_G(\sigma)(\lambda)$ in terms of
the index of an elliptic operator on a reduced space.
However, we believe that our formula provides a unifying framework for all "known cases" of geometric multiplicity formulae and we hope to justify in a future work this statement.

First recall the  definition of a transversally elliptic symbol.
 We denote by $(x,\xi)$, with $x\in M$ and $\xi\in T_x^*M$,
 a point of the cotangent bundle $T^*M$, and by $p: T^*M\to M$ the projection.
 Let $T^*_GM\subset T^*M$  be the set of $(x,\xi)\in T^*M$ such that $\xi$ is orthogonal to the tangent space  to the $G$ orbit through $x$. We say that $\xi$ is a transverse direction.
Let ${\mathcal E}^{\pm}$ be two $G$-equivariant complex vector bundles over $M$. A $G$-equivariant bundle map  $\sigma:p^*{\mathcal
E}^+\to p^*{\mathcal E}^-$  will be called a symbol. The support $\supp(\sigma)\subset T^*M$ of $\sigma$ is the set of elements $(x,\xi)\in
T^*M$ such that the linear map $\sigma(x,\xi):{\mathcal E}^+_{x}\to {\mathcal E}^-_{x}$ is not invertible.
A symbol $\sigma$ is called  elliptic if its support is  compact.  A symbol $\sigma$ is called
transversally elliptic if its support intersected with $T^*_GM$ is  compact.
In other words, $\sigma$ is ``elliptic" in the directions transverse to the orbits.

Let us give some basic examples of transversally elliptic symbols.

$\bullet$ If $D: \Gamma(M,\CE^+)\to \Gamma(M,\CE^-)$ is a $G$-invariant elliptic differential operator on a compact $G$-manifold $M$, acting on the spaces $\Gamma(M,\CE^\pm)$ of smooth sections of $\CE^{\pm}$,
the principal symbol $\sigma$  of $D$ is elliptic, thus transversally elliptic,  and $\Index_G(\sigma)(g)=\Tr_{Ker(D)}(g)-\Tr_{Coker(D)}(g)$.

$\bullet$ If $M$ is homogeneous, any symbol is transversally elliptic.

\bigskip

Let us describe our result. In this introduction, we will make the simplifying hypothesis that all stabilizers of the action of $G$
on $M$ are connected.

If $G$ is a torus with Lie algebra $\g$,
 we parameterize  $\hat G$ by  the lattice $\Lambda\subset \g^*$ of weights of $G$.
 We denote by $g^{\lambda}$ the character of $G$ indexed by $\lambda$.
 Given a transversally elliptic symbol $\sigma$,
 we construct a particular  piecewise polynomial function  $m_G(\sigma)$ on $\g^*$.
 Under the hypothesis that stabilizers of the action of $G$ on $M$ are connected, the function
 $m_G(\sigma)$ is continuous and extends the function
  $\lambda \mapsto \mult_G(\sigma)(\lambda)$ on  $\Lambda$.

A trivial example is when $G$ acts on $G=M$ by left translations and $D$ is the operator $0$, with  index function the trace of
   the regular representation of $G$, that is the $\delta$-function at $1$:
  $$\delta_1(g)=\sum_{\lambda\in
  \Lambda} g^{\lambda}.$$
  Then the function  $\lambda\mapsto \mult_G(D)(\lambda)$ is identically equal to $1$ on $\Lambda$, and is extended by the constant function $m_G(D)=1$.

In the case where $G$ is any compact connected Lie group,
we parameterize $\hat G$ by the set of admissible regular coadjoint orbits (see Section \ref{secnonabelian}).
We construct a $G$-invariant function $m_G(\sigma)$ on $\g^*$,
 and  we determine the multiplicity $\mult_G(\sigma)(\lambda)$ of the irreducible representation $V_\lambda$  parameterized by $G\lambda$ in $\Index_G(\sigma)$ in terms of the value of $m_G(\sigma)$ at $\lambda$.

Our formula for $m_G(\sigma)$ involves  equivariant cohomology  classes on the $G$-manifold  $T^*M$, namely the equivariant Chern character of $\sigma$, and the equivariant truncated $\hat A$  class.
Let $N$ be a $G$-manifold and let $\A(N)$ be the space of differential forms on $N$, graded by its exterior degree.
Following \cite{BV82} and \cite{Wittenequi}, an equivariant form is a $G$-invariant smooth function
$\alpha: \g\to \A(N),$ thus
$\alpha(X)$ is a differential form on $N$ depending differentiably of $X\in \g$.
Consider the operator
\begin{equation}\label{DX}
d_\g\alpha(X)=d\alpha(X)-\iota(v_X)\alpha(X)
\end{equation}
 where $\iota(v_X)$ is the contraction by the vector field $v_X$
generated by the action of $-X$ on $N$.
Then $d_\g$ is an odd operator with square $0$, and  the equivariant cohomology $\CH^*_G(N)$ is defined to be the cohomology space of $d_\g$.
It is important to note that the dependance of $\alpha$ on $X$ may be $C^\infty$.
If the dependance of $\alpha$ in $X$ is polynomial, we denote by
$H^*_G(N)$ the corresponding  $\Z$-graded algebra.
By definition, the grading of $P(X)\otimes \mu$, $P$ a homogeneous polynomial and $\mu$ a differential form on $N$,
is the exterior degree of $\mu$ plus twice the polynomial degree in $X$.

When $\sigma$ is elliptic, $\sigma$ determines an element of the topological equivariant $K$-group $K^0_G(T^*M)$, and its equivariant Chern character
$\Ch(\sigma)(X)$ is a compactly supported
equivariant cohomology class (with $C^{\infty}$ coefficients).
When $\sigma$ is only transversally elliptic, then $\sigma$ determines  an element of $K^0_G(T^*_GM)$
and we have described in \cite{dpv3} a Cartan model for the equivariant cohomology with compact support of $T^*_GM$.
Thus a representative of $\Ch(\sigma)(X)$  is a compactly supported
equivariant form on $T^*M$ (depending on $X$ in a $C^{\infty}$ way), and equivariantly closed on a small neighborhood of $T^*_GM$ in $T^*M$.

Let $J(A)=\det_{\R^d} \frac{e^{A}-1}{A}$, an invariant  function of  $A\in \End(\R^d)$.
 Then, $J(0)=1$. Consider $\frac{1}{J(A)}$ and its Taylor expansion at $0$:
 $$\frac{1}{J(A)}=\det_{\R^d}(\frac{A}{e^{A}-1})= \sum_{k=0}^{\infty}  B_{k}(A).$$
 Each  function $B_{k}(A)$ is an invariant polynomial of degree $k$ on $\End(\R^d)$ and by the Chern Weil construction,
  $B_{k}$
determines an equivariant characteristic class $B_{k}(M)(X)$ on $M$ of homogeneous degree $2k$.

 We define the formal series of equivariant cohomology classes:
   $$[B(M)](X)=\sum_{k=0}^{\infty}  B_{k}(M)(X).$$
   For $X$ small enough, and $M$ compact, the series is convergent, and its sum $B(M)(X)$ is the equivariant  $\hat A$   class of  $T^*M$, that is the square of the equivariant $\hat A$ class of $M$.
   Choose an integer $K$, and define
   $$B(M,K)(X)=\sum_{k=0}^{K}  B_{k}(M)(X).$$
Thus $B(M,K)$ defines an equivariant cohomology class {\bf with polynomial coefficients} on $M$.
We call it the truncated $\hat A$  class of $T^*M$, and $K$ is the order of truncation.

    Let $j_\g(X)=\det_\g \frac{e^{ad X}-1}{ad X}$ and $j_\g^{1/2}(X)$ its square root, a $G$-invariant analytic function on $\g$.

Assume first that $\sigma$ is elliptic.
We define  the $C^{\infty}$ function $I_G(\sigma,K)(X)$ on $\g$ by
$$I_G(\sigma,K)(X)=j_\g^{1/2}(X) \left(\frac{1}{(2i\pi)^{\dim M}}\int_{T^*M}\Ch(\sigma)(X) B(M,K)(X)\right).$$
Consider  its Fourier transform $S_G(\sigma,K)$, a generalized function of $y\in \g^*$.
 In the sense of generalized functions, $S_G(\sigma,K)(y)$ is given by the double integral
 $$S_G(\sigma,K)(y)$$
 $$=
 \frac{1}{(2i\pi)^{\dim M}}\int_{ T^*M}\int_\g j_\g^{1/2}(X)\Ch(\sigma)(X) B(M,K)(X)e^{-i\la y,X\ra} dX.$$

When $\sigma$ is only transversally elliptic,
introduce the canonical Liouville one form  $\ell$ on $T^*M$. Let $\omega=-\ell$.
The  infinitesimal index (as defined in \cite{dpv3}) allows us to  still define the double integral $S_G(\sigma,K)$ above by multiplying
 the  integrand by the oscillatory factor $e^{isd_\g\omega(X)}$ (which is congruent to $1$ in equivariant cohomology). That is:
$$S_G(\sigma,K)(y)$$
 $$=\lim_{s\mapsto \infty} \frac{1}{(2i\pi)^{\dim M}}\int_{ T^*M}\int_\g  e^{is d_\g\omega(X)}j_\g^{1/2}(X)\Ch(\sigma)(X) B(M,K)(X)e^{-i\la y,X\ra} dX.$$
We denote the corresponding limit of the double integral by
$$\frac{1}{(2i\pi)^{\dim M}}\int_{ T^*M}^{\omega}\int_\g  j_\g^{1/2}(X)\Ch(\sigma)(X) B(M,K)(X)e^{-i\la y,X\ra} dX.$$

Assume first that $G$ is a torus. In this case, the function $j_\g$ is identically equal to $1$.
We prove :

$\bullet$ The generalized function
$S_G(\sigma,K)$ has singularities on a union $\CH$ of a finite number of  affine hyperplanes.

$\bullet$ The restriction of $S_G(\sigma,K)$ to each connected component $\c$ of the complement of $\CH$ is a polynomial function.
 This polynomial function is independent of  $K$ if $K\geq
 \dim(M)$.

This allows us to define, for $y\in \g^*\setminus \CH$,
$$m_G(\sigma)(y)=$$
 $$=
 \frac{1}{(2i\pi)^{\dim M}}\int_{ T^*M}^{\omega}\int_\g \Ch(\sigma)(X) \sum_{k=0}^{\infty}B_k(M)(X)e^{-i\la y,X\ra} dX$$
$$=
 \frac{1}{(2i\pi)^{\dim M}}\int_{ T^*M}^{\omega}\int_\g \Ch(\sigma)(X)[B(M)](X)e^{-i\la y,X\ra} dX$$
as the terms associated to $k>\dim M$ vanishes.

The function $m_G(\sigma)(y)$ is a piecewise polynomial function on $\g^*$.
In the language of approximation theory \cite{deboor}, a spline function is a piecewise polynomial function satisfying some further continuity properties.

\begin{theorem}
The function  $m_G(\sigma)$ defined on $\g^*\setminus \CH$ extends to a continuous function on $\g^*$, still denoted
by $m_G(\sigma)$.

If $\lambda$ belongs to the lattice $\Lambda$ of weights of $G$,
we  have the equality
$$\mult_G(\sigma)(\lambda)=m_G(\sigma)(\lambda).$$
\end{theorem}

Thus we have constructed a canonical spline function on $\g^*$, extending the multiplicity function
$\lambda \mapsto \mult_G(\sigma)(\lambda)$.

Let $G$ be a connected compact Lie group.
Let $G\lambda\subset \g^*$ be a regular coadjoint orbit and let ${\rm vol}(G\lambda)$ its symplectic volume.
Similarly, the value of
$S_G(\sigma,K)(y)$ when $y$ is generic and tends to $\lambda$ is well defined, independent of  $K$, provided  $K\geq \dim(M)$,
and
$$\mult_G(\sigma)(\lambda)={\rm vol}(G\lambda)\lim_{y\mapsto \lambda}S_G(\sigma,K)(y).$$

In the general case where stabilizers are not necessarily connected, we define similar
  functions   associated to the action of $G^g$ on $M^g$, where $g$ varies over a finite number of special elements of $G$. Here   $M^g$ is the fixed point submanifold for the action of  $g\in G$ on $M$ and $G^g$ the centralizer of $g$.
 The multiplicity  is described with the help of these functions as a {\bf  piecewise quasi polynomial function}.

\bigskip

We now explain the motivation of this formula.
Let $M$ be a compact $G$-manifold, where $G$ is a torus.
Consider $D$  an elliptic $G$-invariant differential operator, with principal symbol $\sigma$.
 The function $\Index_G(D)(g)=\sum_{\lambda\in \hat G} \mult_G(\sigma)(\lambda)\chi_\lambda(g)$
 is then an analytic function on $G$.
  Recall the formula: for $X\in \g$ small,

  \begin{equation}\label{BV}
  \Index_G(D)(\exp X)=\frac{1}{(2i\pi)^{\dim M}}\int_{T^*M} \Ch(\sigma)(X) B(M)(X),
  \end{equation}
  obtained in (\cite{BV1}, \cite{BV2}),  a "delocalized" version of Atiyah-Bott-Segal-Singer \cite{Atiyah-Segal68}
  equivariant index formula.
  The function $X\mapsto B(M)(X)$ is  defined only when $X$ is sufficiently small.
Inspired by the inversion formulae of box splines \cite{deboor},
 we have replaced the equivariant class $X\mapsto B(M)(X)$
  by its approximation  $B(M,K)(X)$, which is a polynomial function of $X$.
  Under the hypothesis that stabilizers are connected, when $K$ is sufficiently large,  the value of
the Fourier transform of  $\frac{1}{(2i\pi)^{\dim M}}\int_{T^*M} \Ch(\sigma)(X) B(M,K)(X)$ at a integral point $\lambda\in \Lambda$ coincide with the
Fourier coefficient $\mult_G(\sigma)(\lambda)$ of  the periodic function $\Index_G(D)(g)$.

Let us compare the  formulae for $\mult_G(\sigma)$, a function on $\Lambda\subset \g^*$, and  for $m_G(\sigma)$, a function on $\g^*$.
 The number $\mult_G(\sigma)(\lambda)$ is the Fourier coefficient  of  the periodic function $\Index_G(D)(g)$, thus  is given by
$$\mult_G(\sigma)(\lambda)=\int_G \Index_G(D)(g)g^{-\lambda}dg.$$
Replace formally
$\Index_G(D)(\exp X)$ by $\frac{1}{(2i\pi)^{\dim M}}\int_{T^*M}\Ch(\sigma)(X) [B(M)](X)$
and integrate on $\g$ instead of $G$. This gives the similar double integral formula:
$$m_G(\sigma)(y)=\frac{1}{(2i\pi)^{\dim M}}\int_{T^*M}\int_\g \Ch(\sigma)(X) [B(M)](X) e^{-i\la y,X\ra} dX.$$
We have  the miraculous equality

\begin{equation}\label{miraculous}
\mult_G(\sigma)(\lambda)=m_G(\sigma)(\lambda)
\end{equation}
at all integral points $\lambda\in \Lambda.$

This double integral formula for $m_G(\sigma)(\lambda)$ is
reminiscent of Witten non abelian localization formula for computing characteristic numbers of symplectic reductions.
In particular,  Jeffrey-Kirwan \cite{jk} used the truncated Todd class of $M$ in order to compute Riemann-Roch numbers on symplectic reductions.
As pointed out by Guillemin \cite{Guiformal}, the Witten non abelian localization formula  gives an heuristic proof of the $[Q,R]=0$ theorem. This certainly acted as an inspiration for our work. However, our context is different.
%Let us stress again that
% our formula does not provide a proof of the $(Q,R)=0$ theorems of Meinrenken-Sjamaar in the Hamiltonian context, nor of the Paradan-Vergne formula in the spinc context.
 First remark that  we work on $T^*M$  and $\sigma$ is {\bf any}
 transversally elliptic symbol. Our formula for $m_G(\sigma)$ is clearly additive on symbols.
 For general symbols,  we do not have a geometric interpretation of the double integral defining  $m_G(\sigma)(y)$.
  So further work  have  to be done  for reinterpreting $m_G(\sigma)(\lambda)$ in terms of an index on a reduced space,
  for the special symbols of Dirac operators twisted by line bundles, as in \cite{Meinrenken-Sjamaar}, \cite{pep-vergne:spinc}.
 Clearly the symbols have to be very special to obtain such a geometric interpretation. Think to the case
 where  $D$ is the twisted Dolbeaut Dirac operator by a sum of two  lines bundles $L_1,L_2$,  then the multiplicity function $m$ is the sum of the multiplicity functions $m_1,m_2$. Each one of the functions $m_1,m_2$ has a geometric interpretation in terms of reduced spaces, but the geometric meaning for the sum $m$ is not clear.

\bigskip

Let us now give some indications of the proofs.
Following the method of Atiyah-Singer as described in the monograph \cite{At}, we  reduce to the case of the maximal torus $T$ of $G$.
In particular, if $\sigma$ is a $T$-transversally elliptic symbol, we can define the $T$-multiplicity function
 $m_T(\sigma)(y)$, and our formula for the generalized function  $m_G(\sigma)$ is nothing else that the usual formula
$m_G(\sigma)(y)=\sum_w \epsilon(w) m_T(\sigma)(y+w\rho)$, an alternate sum over the Weyl group.
When $M$ is a vector space with a linear action of a torus $G$, the construction of $m_G(\sigma)$ was the object of the articles \cite{dpv1},
\cite{dpv2}. Our proof is very similar. We use Atiyah-Singer description of a set of generators of $K_G^0(T^*_GM)$.
This reduces (almost) our study to the case where $M=P\times V$ where $P$ is a space with a free action of $G$, and $V$ is a linear representation of $G$.  Then the miraculous equality (\ref{miraculous})  follows from a generalized Dahmen-Micchelli inversion formula for box splines (see \cite{vergnepoisson}), a Riemann-Roch formula in approximation theory.
Simple examples of the inversion formula for box splines are given in Example \ref{box4}
and Subsection \ref{amusing}.

\bigskip

\bigskip

It may be useful  to compute our formula  for the spline function
$m_G(\sigma)(\xi)$ on the very simple example of the projective space.

\begin{example} The projective space\end{example}
Let $M=P_1(\C)$ with action of $S^1$ acting on homogeneous coordinates $[z_1,z_2]$ by $[tz_1,z_2]$.
Let $E$ be the element of $Lie(S^1)$ such that $\exp(x E)=e^{ix}\in S^1$.
We identify the lattice of characters of $S^1$ to $\Z$.

Let $0\to L\to \C^2\to Q\to 0$  be the exact sequence of vector bundles on $M$, where $L$ is the tautological bundle and $Q$ its quotient bundle. Let $ a,b$ be integers, and  consider the Dolbeaut  operator $\overline \partial_{a,b}$ acting on sections of the graded vector bundle
$L^a\otimes Q^{b}\otimes \Lambda^{\bullet}\Omega^{0,1}(M)$.
Let $\sigma_{a,b}$ be the principal symbol of $\overline \partial_{a,b}$.
As shown by Atiyah-Bott, the index  of $\overline \partial_{a,b}$ depends only of the two line bundles
$\CL^+=L^a\otimes Q^{b}$  and $\CL^-=L^a\otimes Q^{b}\otimes \Omega^{0,1}(M)$.
The supertrace of the action of $\exp(xE)$ on the fiber of the bundle $\CL^+\oplus \CL^-$
at the fixed point $[1,0]$ is $e^{iax}(1-e^{-ix})$ while it is $e^{ibx}(1-e^{ix})$ at the fixed  point $[0,1]$.
Atiyah-Bott Lefschetz  formula  is
\begin{equation}\label{indexab}
\Index_{S^1}(\overline \partial_{a,b})(e^{ix})=\frac{e^{iax}(1-e^{-ix})}{(1-e^{ix})(1-e^{-ix})}+
\frac{e^{ibx}(1-e^{ix})}{(1-e^{ix})(1-e^{-ix})}.\end{equation}
If $a\leq b$, we obtain
 $$\Index_{S^1}(\overline \partial_{a,b})(e^{ix})=\sum_{k=a}^b e^{ik x}.$$
So the multiplicity function $k \mapsto \mult_{S^1}(\sigma_{a,b})(k)$ on $\Z$ is
$\mult_{S^1}(\sigma_{a,b})(k)=1$ if $0\leq a\leq k\leq b$, or $0$ otherwise.

 Let us now compute our formula.
 Here the equivariant class $[B(M)]$ is the Taylor series of
$\frac{x^2}{(1-e^{ix})(1-e^{-ix})}$:

$$[B(M)](x)=\sum_{k=0}^{\infty} c_{k} x^{k}=1+\frac{1}{12}x^2+\frac{1}{240}x^4+\cdots.$$
%$$(add-bernoulli(2*k,0)*(2*k-1)*(-1)^k*x^(2*k)/(2*k)!)$$

Thus $$I_{S^1}(\sigma_{a,b},K)(xE)=(\sum_{k=0}^{K} c_k x^k)\left(\frac{1}{2i\pi}\int_{T^*M}\Ch(\sigma_{a,b})(xE)\right).$$
 Using the localization formula in compactly supported equivariant cohomology on $T^*M$, we obtain
 $$\frac{1}{2i\pi}\int_{T^*M}\Ch(\sigma_{a,b})(xE)=\frac{e^{iax}(1-e^{-ix})}{x^2}+
\frac{e^{ibx}(1-e^{ix})}{x^2}.$$

By Fourier transform, we obtain
 $$\frac{1}{2i\pi}\int_{T^*M}\Ch(\sigma_{a,b})(xE)=\int_{\R} e^{-i yx} s_{a,b}(y)dy.$$

Here  $s_{a,b}(y)$  is the continuous spline function on $\R$ represented (for $a=0, b=3$) by the  graph of Figure \ref{BoxP1}.
\begin{figure}
\includegraphics[width=3.5 in]{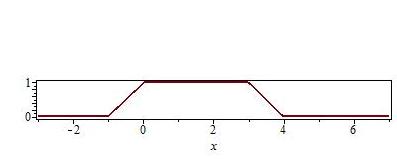}
\end{figure}\label{BoxP1}

Thus the Fourier transform
$S_{S^1}(\sigma_{a,b},K)$  of $I_{{S^1}}(\sigma_{a,b},K)$ is easy to compute:
$$S_{S^1}(\sigma_{a,b},K)(y)=s_{a,b}(y)+ \frac{1}{12} (\delta_a(y)-\delta_{a-1}(y)
+\delta_b(y)-\delta_{b+1}(y))$$
$$- \frac{1}{240} \partial_y^2\cdot(\delta_a(y)-\delta_{a-1}(y)
+\delta_b(y)-\delta_{b+1}(y))+\cdots.$$
We see that
$S_{S^1}(\sigma_{a,b},K)$ coincide with the continuous spline function $s_{a,b}$ on
$\R\setminus\{a-1,a,b,b+1\}$ as soon as $K\geq 1$.
We thus have

$$m_{S^1}(\sigma_{a,b})=s_{a,b}.$$

Remark that  the value of $s_{a,b}$ at an integer $k$ is $0$, except if $a\leq k\leq b$, where the value is $1$.
Our theorem is verified:
$$\mult_{S^1}(\sigma)(k)=m_{S^1}(\sigma_{a,b})(k).$$

\bigskip

Consider now the group $SU(2)$ and let $Z$ be its center.
Let $G=SU(2)/Z$. Then $G$ acts on $M=P_1(\C)$ with connected stabilizers.
We take again the $G$-equivariant line bundle $\CL_{a,b}$  with  $a\leq b$ of same parity (in order that the corresponding representation of $G$
 factorize by $Z$).
 Then  the index of $\overline\partial_{a,b}$ operator is the irreducible representation of $G$, of odd dimension $(b-a)+1$.

Let us now compute our formula.
 Let $\g$ be the Lie algebra of $SU(2)$.
  We identify both $\g$ and $\g^*$ with $\R^3$, with
   $$\g=\{\left(
           \begin{array}{cc}
             ix_1 & x_2+ix_3  \\
             -x_2+ix_3 & -ix_1 \\
           \end{array}
         \right)\}, \hspace{1cm}\g^*=\{\left(
           \begin{array}{cc}
             iy_1 & y_2+iy_3  \\
             -y_2+iy_3 & -iy_1 \\
           \end{array}
         \right)\}.$$

The standard scalar product is $G$-invariant.
A coadjoint orbit is the sphere $S_r= \{y_1^2+y_2^2+y_3^2=r^2\}$.
The Liouville measure on $S_r$ is   such that ${\rm vol}(S_r)=r$, and the admissible coadjoint orbits for $G$ are the spheres of odd radius.
The irreducible representation of $G$ of dimension $k$ is parameterized by the coadjoint orbit $S_{k}$.

Let $$H_\alpha=\left(
           \begin{array}{cc}
             i & 0  \\
             0& -i
           \end{array}
         \right).$$
         Then $\exp(xH_\alpha)$ transforms $[z_1,z_2]\in P_1(\C)$ in $[e^{2ix}z_1,z_2]$.

Consider
$$I_G(\sigma_{a,b},K)(X)=j_\g^{1/2}(X)\int_{T^*M} \Ch(\sigma_{a,b})(X) B(M,K)(X),$$
a $G$-invariant function on $\g$. Let $R=b-a$.
Following the same computation as in the preceding case, we obtain
$$I_G(\sigma_{a,b},K)(xH_\alpha)=\left(\frac{e^{ix}-e^{-ix}}{2ix}\right) V(R,K)(x)$$
with
$$V(R,K)(x)=(\sum_{k=0}^{K} c_k (2x)^{k})\frac{e^{i Rx}(1-e^{2ix})+e^{-i Rx}(1-e^{-2ix})}{4x^2}.$$

Let $r=(y_1^2+y_2^2+y_3^2)^{1/2}$ and let
 $v$ be the $G$-invariant function on $\g^*$ such that
$$v(y_1,y_2,y_3)=\frac{r-(R-1)}{2r},\hspace{1cm}  {\rm if} \,\,R-1\leq r\leq R+1,$$

 $$v(y_1,y_2,y_3)=\frac{(R+3)-r}{2r},\hspace{1cm}  {\rm if} \,\,R+1\leq r\leq R+3$$
 and $0$ otherwise. We see that $v$ is a continuous function on $\g^*$.

  The Fourier transform $S_G(\sigma_{a,b},K)$ of $I_G(\sigma_{a,b},K)$
  coincide with $v$ outside spheres of radius $R-1,R+1,R+3$, as soon as $K>1$.
We thus have $m_G(\sigma)=v$.

Recall that $R=b-a$ is even.
 We see that $v$ vanishes on all spheres of odd radius, except on the sphere of radius $R+1$, parameterizing our representation
 $\Index_G(\overline\partial_{a,b})$, where its value is $1/(R+1)$.
 Thus our formula is verified:

 $$\mult_G(\lambda)={\rm vol}(G\lambda)v(\lambda).$$

\section{Preliminaries}

\subsection{Formal series}\label{defi:JD}
Let $E$ be a vector space.
We introduce $E[[q]]$ as the space of formal series $f[q]=\sum_{k=0}^{\infty} q^k f_k$ with $f_k \in E$.
If the series is convergent ($E$ being a topological vector space) at $q=1$, we write
$f[1]$ for  the sum $\sum_{k=0}^{\infty} f_k$.
If $E$ is an algebra and $f_0$ invertible, we can define $1/f[q]$ in $E[[q]]$.

If $\theta(q,X)$ is a smooth function of $q\in \C$ defined in a neighborhood of $q=0$ and depending  of  a parameter $X$, we denote by
$$\theta([q])(X)=\sum_{k=0}^{\infty} q^k \theta_k(X)$$ its Taylor series at $q=0$. This is a formal series of functions of $X$.

Let $V$ be a real  vector space.
 Let $A\in \End(V)$. Consider the function
  $$J_V(A)=\det_V(\frac{e^{A}-1}{A}).$$
  If $V$ is a complex vector space, we define
   $$J_V^{\C}(A)=\det_{V}^{\C}(\frac{e^{A}-1}{A})$$
   where here the determinant is the complex determinant.
   If $V$ is an Euclidean vector space, and $A$ antisymmetric,
   one has $$J_V(A)=J_{V\otimes_\R \C}^{\C}(A) J_{V\otimes_\R \C }^{\C}(-A).$$
Introduce a variable $q$, and consider
$J_V(q,A)=J_V(q A)$.
Then $$J_V([q])(A)=\sum_{k=0}^\infty q^k T_k(A)=1+q  \frac{\Tr(A)}{2}+\cdots$$
  where $T_k$ is an invariant homogeneous polynomial of degree $k$  on $\End(V)$.
We may also write $J_V([q]A)$ instead of $J_V([q])(A)$.

  If $N$ is a real vector space and $s\in \End(N)$ is a transformation of $N$, we denote by
$GL(s)$ the group of invertible linear transformations of $N$ commuting with $s$. We consider
$$D_N(q,s)(A)=\det_N(1-se^{qA}),$$
an analytic function on $\End(N)$. Write the Taylor series $$D_N([q],s)(A)=\sum_{k=0}^{\infty} q^k D_k^s(A).$$
Then $D_k^s(A)$ is a homogeneous polynomial of degree $k$ on $\End(N)$, invariant by conjugation under $GL(s)$.
If $1-s$ is invertible, we can define $1/D_N([q],s)(A)$.

\subsection{Transversally elliptic symbols}
Let $G$ be a compact Lie group. Let $\g$ be its Lie algebra, and $\g^*$ be the dual vector space of $\g$.
If $G$ acts on a $G$-manifold $N$, and $X\in \g$,  we denote by $$v_X(n)=\frac{d}{d\epsilon} \exp(-\epsilon X)\cdot n|_{\epsilon=0}$$ the vector field on $N$ generated by $-X$.

 Let us consider two $G$-equivariant hermitian bundles $\CF^\pm$ on $N$ and
let $\sigma:\CF^+\to \CF^-$ be a $G$-equivariant morphism. The support $\supp(\sigma)$  of $\sigma$ is the set of points $n\in N$ where
$\sigma_n$ is not invertible. We will say that $\sigma$ is elliptic if $\supp(\sigma)$ is compact.  Then $\sigma$ determines an element of the $G$-equivariant topological $K$-group
$K_G^0(N)$, still denoted by $\sigma$.

Recall the definition of the Bott symbol.
 Let $N$ be a   Hermitian  vector space with Hermitian form  $\la\cdot,\cdot \ra$ and complex structure $J$.
  We denote by $\mathrm{U}$ the unitary group of transformations of $N$, and by $U(1)\subset \mathrm{U}$ the subgroup  formed the homotheties $n\mapsto e^{i\theta} n$.   Consider $S=\bigwedge N$, graded in even and odd degree. Let ${\bf c}:N \to  \End(S)$
  be the Clifford action. The map $\bf c$
   is equivariant, interchanges $S^\pm$ and   ${\bf c}(n)^2=-\|n\|^2 {\rm Id_S}.$
    If $N=\C$, then
  $${\bf c}(z)=\left(
               \begin{array}{cc}
                 0 & -\overline{z} \\
                 z & 0 \\
               \end{array}
             \right).$$

\begin{definition}\label{defi:Bott}
The Bott symbol $Bott(N,J) \in \Gamma(N, hom(S^+,S^-))$  is defined by
$$
Bott(N,J)(z)={\bf c}(z): \bigwedge^{even}_J N\longrightarrow \bigwedge^{odd}_JN, \qquad z\in N.
$$
\end{definition}
It is an elliptic symbol, with support $\{0\}$, and equivariant with respect to $\mathrm{U}$.

If we take the opposite complex structure on $N$, one has the relation
$Bott(N,-J)=(-1)^{\dim_\C N}\chi Bott(N,J)$ in $K_{\mathrm{U}}^0(N)$, where $\chi$ is the character $\chi(g)=(\det^{\C}_N(g))^{-1}$
of $\mathrm{U}$. This is easily seen by restriction to $\{0\}\mapsto N$.

%More generally, if $\CN\to M$ is a Hermitian vector bundle with complex structure $J$,
%we can define the Bott symbol $$Bott(\CN,J)(x,z)
%={\bf c}(z): \bigwedge^{even} \CN_x\longrightarrow \bigwedge^{odd}\CN_x, \qquad x\in M, z\in \CN_x.
%$$

\bigskip

Let $M$  be a $G$-manifold,  $N=T^*M$  its cotangent bundle.
 We denote by $(x,\xi)$, with $x\in M$ and $\xi\in T_x^*M$,
a point of the cotangent bundle $T^*M$, and by $p: T^*M\to M$ the projection.
Let $\ell$ be the Liouville one form: for $x\in M$, $\xi\in T_x^*M$ and $V$ a tangent vector at the point $(x,\xi)\in T^*M$, $
\ell_{x,\xi}(V)=\ll \xi,p_*V\rr.$ Then $\Omega=-d\ell$ is the symplectic  form of $T^*M$, and we  use
the corresponding orientation of $T^*M$ to compute integrals
  on $T^*M$
  of differential forms with compact support.
Denote by $T^*_GM\subset T^*M$ the union  of  the space of covectors conormal to the $G$-orbits. Let $\mu: T^*M\to \mathfrak g^*$ be the moment map  $\ll \mu(x,\xi),X\rr=\la \xi, v_X\ra$.
Then $T^*_GM=\mu^{-1}(0)$ is the zero fiber  of the moment map $\mu$.

Let ${\mathcal E}^{\pm}$ be two $G$-equivariant complex vector bundles over $M$. So a $G$-equivariant bundle map  $\sigma:p^*{\mathcal
E}^+\to p^*{\mathcal E}^-$ is a symbol.   A symbol $\sigma$ is called
transversally elliptic if its support intersected with $T^*_GM$ is  compact.
If $\sigma$ is transversally elliptic, it determines an element of $K_G^0(T^*_GM)$, still denoted by $\sigma$.
Atiyah-Singer (see the monograph \cite{At}) have associated to any
element $\sigma\in K_G^0(T^*_GM)$ a   virtual trace class representation of   $G$.
 Let
$\Index_G(\sigma)(g)$ its trace:
$$\Index_G(\sigma)(g)=\sum_{\lambda\in \hat G} m_G(\sigma)(\lambda) \chi_\lambda^G(g).$$
Here $\chi_\lambda^G$  is the trace of the unitary irreducible representation $V_\lambda^G$ of $G$.
When $G$ is given, we might write simply $\chi_\lambda$ or $V_\lambda$ instead of
$\chi_\lambda^G$ or $V_\lambda^G$.
We might also write
$$\Index_G(\sigma)=\oplus_{\lambda\in \hat G} m_G(\sigma)(\lambda) V_\lambda^G,$$
an infinite sum of irreducible representations with finite multiplicities
$m_G(\sigma)(\lambda)\in \Z$.

Recall examples of transversally elliptic symbols.

\begin{example}\label{ex:elliptic} Elliptic symbols \end{example}

Any elliptic symbol is transversally elliptic.
 If $M$ is an even dimensional compact manifold and is oriented, then any element of $K_G^0(T^*M)$ is the symbol of a twisted Dirac operator.

\begin{example}\label{ex:homogeneou} Branching rule\end{example}

Let $M=G$ and let $H\subset G$ be a compact subgroup of $G$. Consider the action of $G\times H$   on $M$  by left and right translations.
Let $\CE^+=M\times \C$ be the trivial vector bundle, and $\CE^-=M\times\{0\}$.
 Then the $0$ symbol $\sigma_0$ is transversally elliptic with respect to the action of $G\times H$ on $M$, and
$$\Index_{G\times H}(\sigma_0)(g,h)=\sum_{\lambda\in \hat G,\mu\in \hat H} m_{G\times H}(\sigma_0)(\lambda,\mu) \overline{\chi_\lambda^G(g)}\chi_\mu^H(h),$$ where
$m_{G,H}(\sigma_0)(\lambda,\mu)$ is the multiplicity of the irreducible representation
$V_\mu^H$ of $H$ in the irreducible representation $V_\lambda^G$ of $G$.
 The function $m_{G,H}(\sigma_0)(\lambda,\mu)$ is thus the branching function.

\begin{example} Orbifold index\label{ex:orbifold}\end{example}
Let $M$ be a compact manifold  and let $G$ be a compact
group acting on $M$. Assume that all stabilizers $G_m$ of points of $M$ are finite. Then $M/G$ is an orbifold.
Let $\sigma$ be a $G$-transversally elliptic symbol on $M$, and write
$\Index_G(\sigma)(g)=\sum_{\lambda\in \hat G} m_G(\sigma)(\lambda) \chi_\lambda(g)$.
If $\lambda_0$  is the trivial representation of $G$,
 $m_G(\sigma)(\lambda_0)$
is the index of the  elliptic symbol on the orbifold $M/G$ associated to $\sigma$.

\begin{example}  Atiyah symbol \label{pushed symbol} \end{example}
Let $N$ be an even dimensional  vector space
with a linear action of a torus $G$.
We choose a $G$-invariant Euclidean product $\la \cdot,\cdot \ra$.
This allows us to identify $T^*N$ with $N\oplus N$.
Let us choose a $G$ invariant complex structure $J$ on $N$ preserving the inner product and let us consider ${\rm Sym}(N)=\oplus_{k=0}^{\infty} {\rm Sym}^k(N)$
where ${\rm Sym}^k(N)$ is the subspace of $\otimes^k N$ formed of symmetric tensors.
Here $N$ is considered as a complex vector space via the complex structure $J$.
We can consider the Bott symbol $Bott(N,-J)$ associated to the opposite complex structure.
Consider the equivariant map
$q:N\oplus N\to N$ defined by
$(x,\xi)\mapsto \xi+J x$.
\begin{definition}\label{defiat}
The Atiyah symbol $at(N)$  is the reciproc image  of $Bott(N,-J)$ by the map $q$.
\end{definition}
The following proposition is proved in \cite{At} (see also \cite{BV1}).
\begin{proposition}
 $\Index_{\mathrm{U}}^N(at(N))={\rm Sym}(N)$.
\end{proposition}

Similarly, if $\CN\to M$ is a Hermitian vector bundle on $M$,  using a connection, we can define a symbol
$at(\CN)$ on $T^*\CN$.

\subsection{Equivariant cohomology and infinitesimal index}

We recall the definitions of  the equivariant cohomology given in the introduction.
Let $N$ be a $G$-manifold and let $\A(N)$ be the space of differential forms on $N$, graded by its exterior degree.
Consider the operator (Equation (\ref{DX}))
$$d_\g:C^{\infty}(\g,\A(N))^G\to  C^{\infty}(\g,\A(N))^G.$$
The equivariant cohomology $\CH^*_G(N)$ is the cohomology space of $d_\g$.

If the dependance of $\alpha$ in $X$ is polynomial, this model for equivariant cohoùmology
 is equivalent to the Cartan model for topological equivariant cohomology of the $G$-space $N$ (see \cite{cartan}). We denote by
$H^*_G(N)$ the corresponding  $\Z$-graded algebra.

The equivariant integration $\int_N \alpha(X)$ associates to an equivariant cohomology class with compact support on $N$ a
$G$-invariant  $C^{\infty}$ function on $\g$.

If $\CN\to M$ is an oriented $G$-equivariant vector bundle over $M$, the equivariant Thom class  $Thom(\CN)(X)$ of $\CN$ is the equivariant
class with polynomial coefficients, and  compact support along the fiber, such that its integral along the fiber is identically equal to $1$.

In our article, we need to take Fourier transforms of equivariant integrals.
To this purpose,
we  consider  (as in \cite{dpv3})  the subcomplex of    equivariant forms  $\alpha(X)$,  with compact support on $N$,
 which can be written $\alpha(X)=\int_{\g^*} e^{i \la y, X\ra} q(y)$ where $q(y)$ is a distribution with compact support on $\g^*$ with value
 in the space of differential forms with compact support on $N$.
 We denote by $\CH_{G,c}^m(N)$ the corresponding cohomology space (the letter $m$ is for moderate growth).
 The space $\CH_{G,c}^m(N)$ is a module for $H_G^*(N)$.

  Let us consider two $G$-equivariant hermitian bundles $\CF^\pm$ on $N$ and
let $\sigma:\CF^+\to \CF^-$ be a symbol.
We recall the definition of the equivariant Chern character of  $\sigma$.
We choose $G$-invariant hermitian connections $\nabla^{\pm}$ on $\CF^{\pm}$ with curvature $R^\pm$. Then
$R^\pm$ is a  $2$-form on $N$ with values endomorphisms of $\CF^{\pm}$.
Let $\mu^{\CF^{\pm}}(X)$ be the corresponding moment maps determined by the Kostant equation
$$\CL^{\CF^{\pm}}(X)=\nabla^{\pm}_{v_X}+\mu^{\CF^{\pm}}(X)$$
and let $R^\pm(X)= \mu^{\CF^{\pm}}(X)+R^{\pm}$ be the corresponding equivariant curvatures of $\CF^{\pm}$.
We assume that, outside a small neighborhood of the support of $\sigma$, the connections $\nabla^{\pm}$  are transformed to each other by the isomorphism $\sigma$. We say that $\nabla^{\pm}$ are adapted to $\sigma$.
 It follows that
  the closed  equivariant  form $\Ch(\sigma,\nabla)(X):=\Tr (e^{R^+(X)})-\Tr(e^{R^-(X)})$ is a differential form on $N$
supported on a neighborhood of $\supp(\sigma)$.

Let $g$ be an element of $G$, and let $N^g$  be the fixed point submanifold of the action of $g$ on $N$.
Then $g$ acts fiberwise on the bundles $\CF^{\pm}$ restricted to $N^g$.
Let $G^g$ be the centralizer of $g$, and $\g^g$ its Lie algebra.
We restrict the equivariant curvatures $R^\pm(X)$  of $\CF^{\pm}$ to $N^g$.
Then, for $X\in \g^g$, we define
 $\Ch(g,\sigma,\nabla)(X):=\Tr( ge^{R^+(X)})-\Tr(g e^{R^-(X)})$. This is a closed $G^g$-equivariant differential form on $N^g$
supported on a neighborhood of $\supp(\sigma)\cap N^g$.

Assume that $\sigma$ is elliptic, then $\Ch(\sigma,\nabla)(X)$ is a differential form on $N$ with compact support.
 Furthermore, as the dependance of $R^{\pm}(X)$ in $X$ is through the linear map
 $\mu^{\CF^{\pm}}(X)$ with values skew hermitian endomorphisms on $\CF^\pm$, it is easy to see (see Lemma 5.3, \cite{dpv2})
 that we can write  $\Ch(\sigma,\nabla)(X)=\int_{\g^*} e^{i\la y, X\ra} q(y,n)$
 where, for each $n\in N$, $q(y,n)$ is a distribution with compact support on
 $\g^*$. So $\Ch(\sigma,\nabla)$ determines a class $\Ch(\sigma)$ in $\CH_{G,c}^m(N)$, depending only of the class of $\sigma$
   in $K_G^0(N)$.
This is the equivariant Chern character. Similarly  $\Ch(g,\sigma,\nabla)$ determines a class
 $\Ch(g,\sigma)$ in $\CH_{G^g,c}^m(N^g)$, called the twisted equivariant Chern character.

Let us give a simple example (see the  proof of Proposition 5.9 \cite{dpv2}) of the equivariant Chern character of an elliptic symbol.
 Let $N=\C$ be a one dimensional complex vector space, with action of $G=U(1)$. We identify  $\g$ to $\R$ by choosing a basis $E$
 of $\g$ so that $\exp(\theta E)=e^{i\theta}$.
 We consider the Bott symbol $\sigma$ on $N$ (\ref{defi:Bott}).
 Thus $\CF^+=N\times \C$ is the trivial vector bundle on $N$, and  $\CF^-=N\times N$, with morphism
 $\sigma(z): \CF^+_z\to \CF^-_z$ the multiplication by $z$.
 Let us compute its equivariant Chern character.
Let $\chi$ be a function on $\R$ with compact support contained in $|t|\leq 1$ and identically $1$ near $0$.
 Let $\beta=(\chi(|z|^2)-1) \frac{d z}{ z}$.
  This is a well defined one form on $\C$, invariant under the action of $U(1)$.
  The connections $\nabla^+=d$ and $\nabla^-=d+\beta$ are $G$ invariant connections on $\CF^+$, $\CF^-$ and for $|z|>1$, we have
  $\nabla^-\sigma=\sigma\nabla^+$.
  We see that
   $$\Ch(\sigma,\nabla)(\theta )=1-e^{i\theta \chi(|z|^2)}+ e^{i\theta \chi(|z|^2)}\chi'(|z|^2) dz \wedge d\overline z$$
   is a compactly supported equivariant form on $\C$.
  For each $z$, the Fourier transform of $\Ch(\sigma,\nabla)(\theta)$ is supported at the point $-\chi(|z|^2)$.
The integral $\int_\C \Ch(\sigma,\nabla)(\theta)$ is easily computed in polar coordinates.
We obtain
$$\frac{1}{2i\pi}\int_\C \Ch(\sigma,\nabla)(\theta)=\frac{e^{i\theta}-1}{i\theta}.$$
Its Fourier transform is the characteristic function of the interval $[-1,0]$.

We have the following proposition  (Proposition 5.10, \cite{dpv2}).

\begin{proposition}\label{cherntwistedBott}
Let $N$ be a Hermitian vector space, and $G$ be a group acting unitarily on $N$. Then, for any $g\in G$,
we have the equality
$$\Ch(g,Bott(N,J))(X)$$
$$=(2i\pi)^{\dim_{\C}{N^g}}\det_{N^g}^{\C} \left(\frac{e^X-1}{X}\right)\det_{N/N^g}^{\C}(1-g e^X) Thom(N^g)(X)$$
in $\CH_{G,c}^m(N^g)$.
\end{proposition}

\begin{remark}\label{defichBV}
When $\CN$ is  a vector bundle over $M$, and $\sigma$ an elliptic symbol on $\CN$,  we may privilege
a representative of $\Ch(\sigma)(X)$ via
rapidly decreasing differential forms on the fibers of  $\CN$. This is the choice we made in \cite{BV2}.
It is easy to see that the
             construction of equivariant Chern characters via adapted  connections
              differ from BV construction via a boundary which is also rapidly decreasing at $\infty$ on the fibers.
\end{remark}

If $Z$ is a $G$-invariant closed subset  of $N$, we have defined in \cite{dpv3} a Cartan model for the space $\CH_{G}(Z)$ of
equivariant cohomology.
   A representative  is an equivariant form  $\alpha:\g\to \mathcal A(N)$ such that $d_\g\alpha=0$ in a neighborhood of $Z$.
   The dependance of $\alpha$ in $X$ is
   $C^{\infty}$.
   We denote by $1_Z\in \CH_{G}(Z)$ the class represented by a $G$-invariant function
    ${\rm One}_Z$ identically equal to $1$ in a neighborhood of $Z$, and supported on a small neighborhood of $Z$.

We  similarly consider the subcomplex of    equivariant forms  $\alpha(X)$  with {\bf compact support} on $N$
such that $d_\g\alpha=0$ in a neighborhood of $Z$, and
 whose Fourier transform in $X$ is  supported on  a compact subset of $\g^*$, and
we denote the corresponding cohomology group by $\CH_{G,c}^{m}(Z)$.

Let $\sigma$ be a symbol on $N$ and assume that the support of $\sigma$ intersected with $Z$ is compact.
Thus $\sigma$ defines a class  in $K^0_G(Z)$ still denoted by $\sigma$.
We then define a representative of
$\Ch(\sigma)(X)$ to be the equivariant cohomology class ${\rm One}_Z \Ch(\sigma,\nabla)(X)$. Here $\nabla$ and ${\rm One}_Z$ are chosen such that
the support of the differential form ${\rm One}_Z \Ch(\sigma,\nabla)(X)$ is compact.
Similarly, for $g\in G$, we obtain the twisted Chern character $\Ch(g,\sigma)\in \CH_{G^g,c}^{m}(Z^g)$.

We now recall the definition of infinitesimal index \cite{dpv3}.

   Let $\omega$ be a $G$-invariant one form  on $N$, and  $\mu(n)(X)=-\la \omega, v_X\ra$ be the corresponding moment map
   $N\to \g^*$. Let $d_\g\omega(X)=\mu(X)+ d\omega.$
The main example we will consider is when $N=T^*M$, and $\omega=-\ell$ is minus the Liouville form.
So $(d_\g\omega)(X)=\mu(X)+\Omega$   is the equivariant symplectic  form of $T^*M$.

Assume that $Z=\mu^{-1}(0)$.
   Let $\alpha(X)$ be an equivariant form  with compact support on $N$, whose Fourier transform in $X$ is supported on a compact subset $K$ of $\g^*$.
   Then
   $I(s)=\int_{\g} e^{is d_\g\omega(X)} \alpha(X) dX$
   is a differential form on $N$ with support contained on the set of elements $n$ such that $s\mu(n)\in K$.
   So, when $s$ tends to $\infty$, $I(s)$  is supported in a neighborhood of $Z$.
   Consider an element $\theta\in \CH_{G,c}^m(Z)$
   and denote still by $\theta$ a representative of the class $\theta$.
   The infinitesimal index $\Infdex_G^{\omega}(\theta)$ of $\theta$ is the distribution   on $\g^*$
   such that, for a test function $f$ on $\g^*$,
    $$\la \Infdex_G^{\omega}(\theta),f\ra=\lim_{s\mapsto \infty}\int_N \left(\int_\g e^{i s d_\g\omega(X)} \theta(X) {\hat f}(X)  dX\right).$$
This distribution depends only of the class of $\theta$ in $\CH_{G,c}^{m}(Z)$ and not of its representative.
It however depends of $\omega$. We will denote the limit
$\lim_{s\mapsto \infty}\int_N \left(\int_\g e^{i s d_\g\omega(X)} \theta(X) {\hat f}(X)  dX\right)$ by $\int_N^{\omega}\int_\g  \theta(X) {\hat f}(X)dX.$
Thus we write
$$\la\Infdex_G^\omega(\theta),f\ra=\int_N^{\omega}\int_\g  \theta(X) {\hat f}(X)dX.$$

Here ${\hat f}(X)=\int_{\g^*}e^{-i\la y,X\ra} f(y) dy$
and the measure $dX$ is chosen such that $f(y)=\int_\g  e^{i\la y, X\ra}{\hat f}(X)dX.$

Let us give an example of infinitesimal index (see \cite{dpv2}, proof of Theorem 4.21).
Let us consider $N=\R^2=\C$, with action of $G=U(1)$.
Let $Z=\{0\}$, and $\alpha=\frac{1}{2} (x_1 dx_2-x_2dx_1)$ with moment map $\frac{1}{2}(x_1^2+x_2^2)$.
\begin{lemma}\label{exinf} Let $1_Z\in \CH_{G,c}(Z)$ and $f$ be a test function on $\g^*$.
Then
$$\frac{1}{2i\pi}\la\Infdex_G^{\alpha}(1_Z),f\ra=\int_0^{\infty} f(y) dy.$$
\end{lemma}

\subsection{Formal series of equivariant classes}
Let $M$ be our $G$-manifold.

Let $\CV\to M$ be a real or complex vector bundle $G$-equivariant vector bundle on $M$  with typical fiber a real or complex vector space $V$.
 The Chern-Weil map $W$ associates to an $GL(V)$  invariant polynomial $f$ on $\End(V)$ an equivariant characteristic class
 $W(f)$ in $H_G^*(M)$. If $f$ is homogeneous of degree $k$, then $W(f)$ is homogeneous of degree $2k$.
 Our conventions for the Chern-Weil homorphism  $W$   are  as in \cite{B-G-V}.

 Let $A\in \End(V)$. Introduce a variable $q$,  and consider the Taylor expansion
$$J_V([q]A)=\det_V( \frac{e^{[q]A}-1}{[q]A })=\sum_{k=0}^{\infty}q^{k}T_k(A).$$
Our main new concept is the introduction of the following formal equivariant characteristic class of $M$.
\begin{definition} The formal
$J$-class of $M$ is the series of elements of $H_G^*(M)$  defined by
$$J([q],M)=\sum_{k=0}^{\infty} q^{k} W(T_k)$$ obtained by applying
the Chern-Weil map for the real vector bundle $\CV=TM\to M$ to the series
 $\det_V( \frac{e^{[q]A}-1}{[q]A})=\sum_{k=0}^{\infty}q^{k}T_k(A)$.
  \end{definition}
Here $W(T_k)$ is homogeneous of degree $2k$.

We can thus consider  $B([q],M)=1/J([q],M)$ in the ring of formal series of equivariant cohomology classes with polynomial coefficients. We write $B([q],M)=\sum_{k=0}^{\infty} q^kB_k(M)$.
 In the introduction, we have introduced
 $[B(M)]=\sum_{k=0}^{\infty} B_k(M)$ and the truncated class
 $B(M,K)= \sum_{k=0}^{K} B_k(M)$. It is more convenient not to fix an order of truncation, and to work with the full
 formal series $1/J([q],M)$.

When  $G=\{1\}$, then $p^*J([1],M)=J(M)$ is the inverse  of the usual Todd class of the tangent bundle to $T^*M$ (considered as an
almost complex manifold). Furthermore, if $\sigma$ is elliptic, Atiyah-Singer formula \cite{Atiyah-Singer-1} for $\Index(\sigma)\in \Z$ is
$$\Index(\sigma)=\frac{1}{(2i\pi)^{\dim M}}\int_{T^*M}\frac{\Ch(\sigma)}{J(M)}.$$

Consider $$D_N([q],s)(A)=\det_N(1-se^{[q]A})=\sum_{k=0}^{\infty} q^{k} D_k^s(A).$$

Consider the normal bundle $\CN\to M^g$. Thus $g$ produces an invertible linear transformation of $\CN_x$ at any $x\in M^g$. The Chern
Weil homomorphism for the real vector bundle $\CN$ (with structure group $GL(g)$) produces a series $D([q],g,M/M^g):=\sum_{k=0}^{\infty} q^{k}
W(D_k^g)$ of closed equivariant differential forms on $M^g$. The coefficient in $q^0$ of this series is just the function $x\mapsto \det_{\CN_x}(1-g)$, a
function which is a non zero constant on each connected component of $M^g$.

We thus obtain a formal series  of $G^g$ equivariant classes on $M^g$
by considering  the form $J([q],M^g)(X) D([q],g,M/M^g)(X)$.
We can also invert this series in the ring of formal series.

\subsection{Piecewise polynomial functions}\label{subsecPW}
In this subsection, $G$ is  a  torus.
 Let $V=\g^*$  equipped with the lattice $\Lambda\subset \g^*$ of weights of $G$.
If $g=\exp X$, we denote by $g^{\lambda}=e^{i\langle \lambda,X\rangle}$.
The function $g\mapsto g^{\lambda}$ is a character (a one dimensional representation) of $G$.
%We denote by ${\mathcal C}(\Lambda)$ the space of (complex valued) functions on $\Lambda$. %If $g\in G$, we denote by $\hat g$ the
%function $g\to g^{\lambda}$ on $\Lambda$. If $m\in {\mathcal C}(\Lambda)$, then $\hat g m\in \mathcal C(\Lambda)$ is  defined by $(\hat
%g m)(\lambda)= g^{\lambda} m(\lambda)$.

Using the Lebesgue measure $dy$ associated to $\Lambda$, we identify generalized functions on $\g^*$ and distributions  on $\g^*$.
If $h$ is a generalized
function on $\g^*$, we denote by $\int_{\g^*}h(y)f(y)dy$ its value on the test function $f$.
% We denote by $\delta_v$ the Dirac function at the point $v\in V$.

Let  $\CH$  be a finite collection of rational affine hyperplanes in $\g^*$.
An element of $\CH$ will be called an admissible wall. For $H$ an affine hyperplane,
let $lin(H)$ be the hyperplane parallel to $H$.
An element $v\in V$ is called $\CH$-\textit{generic}   if $v$ is not on any hyperplane of the collection $\{lin(H), H\in\CH\}$.  We just say that $v$ is
generic.
A  tope $\c$ is a connected
component of the complement of all admissible walls (thus $\c$ is  an open  convex  subset of $V$) and we denote by $V_{reg}$ the union of topes.
If $v\in V$, and $\epsilon$ is a generic vector, then $v+t\epsilon$ is
in $V_{reg}$ if $t>0$ and sufficiently small.
Assume that $f$ is a function on $V_{reg}$, given on each tope by the restriction of an  analytic function of $y\in V$ (depending of the tope).
We say that $f$ is a piecewise analytic function on $V_{reg}$.

\begin{definition}\label{def:rightlimit} Let $v\in V$, and $f$  a piecewise analytic function. Let $\epsilon$ be a
generic vector. Define $(\lim_{\epsilon}f)(v)=\lim_{t>0,t\mapsto 0}f(v+t\epsilon)$. \end{definition}

A piecewise polynomial function is a function on $V_{reg}$ which is given by a polynomial formula on each tope.
 We denote by $PW$ the
space of piecewise polynomial functions.
 Consider $f\in PW$ (defined on $V_{reg}$) as a locally $L^1$-function on $V$,  thus $f$ defines  a generalized function on $V$.
 An element of $PW$, considered as a generalized function on $V$,
will be called a piecewise polynomial generalized function.
\begin{definition} The space $\CS$ is the space of generalized functions  on
$V$ generated by the action of constant coefficients differential operators  on piecewise polynomial generalized functions.
\end{definition}
 For example, the Heaviside function on $\R$ is a piecewise polynomial generalized function. Its derivative in the sense of generalized functions is
 the Dirac function at $0$, and belongs to $\CS$.
%

%Piecewise analytic $PA$.
%
%A function in $PA$  can be evaluated at any point of $V_{reg}$. If $v\in V$, and $\epsilon$ is a generic vector, then $v+t\epsilon$ is
%in $V_{reg}$ if $t>0$ and sufficiently small.
%

Introduce formal series $m([q])=\sum_{k=0}^{\infty} q^k m_k$ of generalized functions on $V$.
Then, if $f$ is a test function,
$$\int_{\g^*}m([q])(y)f(y)dy=\sum_{k=0}^{\infty} q^k \int_{\g^*} m_k(y)f(y)dy$$ is a formal power series in $q$. It may be
evaluated at $q=1$ if the preceding series is finite (or convergent).
If $\epsilon$ is generic, $v\in V$, and the functions $m_k$ restricted to $V_{reg}$ piecewise analytic, we may define
$$\lim_{\epsilon}m([q])(v)=\sum_{k=0}^{\infty} q^k \lim_{\epsilon} m_k(v).$$
This is a formal series in $\C[[q]]$.

An important piecewise polynomial function is the box spline \cite{deboor}.
Let $\Phi=[\phi_1,\phi_2,\ldots, \phi_N]$ be a list of non zero elements of $\Lambda$.
We assume that the set $\Phi$ spans $V$.
Let $Z(\Phi)=\{\sum_{i=1}^N  t_i \phi_i; 0\leq t_i \leq 1\}$
be the zonotope determined by $\Phi$.
The box spline $B_\Phi(y)$ is the measure on $V$ supported on $Z(\Phi)$ obtained by convolutions of the intervals $[0,1]\phi_k$:
$$\int_V B_\Phi(y) f(y)dy=\int_0^1\cdots \int_0^1 f(\sum_{k=1}^N t_k\phi_k)dt_1\cdots dt_N.$$
The Fourier transform  of $B_\Phi$ is the function
$J_\Phi(X)=\prod_{\phi\in \Phi} \frac{1-e^{-i\la \phi,X\ra }}{i\la \phi,X\ra }$ of $X\in \g$.

Let $\CH_0$ be the set of hyperplanes of $V$ generated by $\dim V-1$ linearly independent elements of $\Phi$.
 Let $F$ be the set of sums $\phi_I=\sum_{i\in I} \phi_i$ of elements of $\Phi$, where $I$ is a sublist of $\Phi$.
 We consider the finite set  $\CH$ of affine hyperplanes of the form $p+H$ with $H\in \CH_0, p\in F$.
 Consider the set $V_{reg}$ of elements $v$ of $V$ such that $v$ does not belong to any affine hyperplane $H$ with $H\in \CH$.
Then $B_\Phi$ is a locally polynomial function with respect to $\CH$.

\begin{example} \label{box4} The box spline for $\Phi=[1,1,-1,-1]$

Consider the function $j_\Phi(x)=(\frac{e^{ix }-1}{ix})^2 (\frac{1-e^{-ix }}{ix})^2$
of $x\in \R$. Its Fourier transform (the convolution of the intervals $[0,1]$, $[0,1]$
$[-1,0]$, $[-1,0]$) is the measure $b_4(y)$
given by
$$
\begin{cases}
   \frac{1}{6}(y+2)^3\hspace{21mm}{\rm if}\ -2<y<-1 ,\\
    \frac{2}{3}-y^2-\frac{1}{2} y^3\hspace{14mm} {\rm if}\ -1\leq y\leq 0,\\
    \frac{2}{3}-y^2+\frac{1}{2} y^3\hspace{14mm} {\rm if}\ 0\leq y\leq 1,\\
    \frac{-1}{6}(y-2)^3\hspace{20mm}{\rm if}\ 1<y<2 .\\
    \end{cases}
$$

Consider the equation
$\frac{1}{j_\Phi(x)} j_\Phi(x)=1$.
Replace $\frac{1}{j_\Phi(x)}$ by its Taylor series
$$\frac{1}{j_\Phi(x)}=1+\frac{1}{6}x^2+\cdots.$$
The  inversion formula for the box spline (due to Schoenberg \cite{sc} for splines in one variable) states that
the Fourier transform $Db_4$
of the function $(1+\frac{1}{6}x^2)j_\Phi(x)$
is continuous and its value at any integer $k$ is equal to $0$, except for $k=0$, where the value is $1$.
The function $b_4$  and the function $Db_4=b_4-\frac{1}{6} (\frac{d}{dy})^2 b_4$
are plotted in the figure below.

\includegraphics[width=6cm]{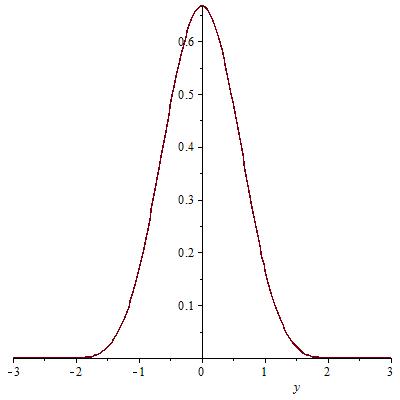}\includegraphics[width=6cm]{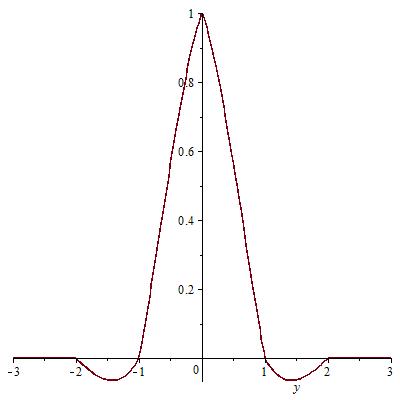}
\label{Box4}.

A more elaborate example is given in Subsection \ref{amusing}.

\end{example}

Assume that the cone $\Cone(\Phi)$ generated by $\Phi$ is a  salient cone in $V$.
Then we can define the spline $T(\Phi)$. This is  is the measure on $V$ supported on $\Cone(\Phi)$ obtained by convolutions of the half lines  $[0,\infty]\phi_k$:
$$\int_V T(\Phi)(y) f(y)dy=\int_0^\infty\cdots \int_0^\infty f(\sum_{k=1}^N t_k\phi_k)dt_1\cdots dt_N.$$
It is a locally polynomial measure.
More generally, if ${\bf y}\in  \C^N$, we introduce the spline function $T(\Phi,{\bf y})$ with parameters ${\bf y}=(y_1,y_2,\ldots, y_N)$
as the convolution of the distributions
 $\int_0^\infty f(t\phi_k)e^{t y_k}dt.$
This function is piecewise analytic.

Let $N$ be a Hermitian vector space  with a linear representation of $G$.
We write $N=\oplus_{\phi\in \Phi} N_\phi$, as a sum of $1$-dimensional representation spaces of $G$ .
If $\Phi$ spans $V=\g^*$, the action of $G$ on $N$ is with finite generic stabilizer.
We still assume that $\Phi$ span a salient cone $\Cone(\Phi)$ in $V$.
Then the action of $G$ in $\Sym(N)$ is with finite multiplicity.
We write $$\Sym(N)=\oplus_{\lambda\in \hat G} \Sym_\lambda(N),$$ its decomposition  with respect to
$G$.
Thus, under the action of $G$, we have
$$\Tr_{\Sym(N)}(g)=\sum_{\lambda \in \hat G} \dim(\Sym_\lambda(N)) g^{\lambda}$$
and the function $\lambda \mapsto \dim(\Sym_\lambda(N))$ is the so called Kostant partition function (for $\Phi$).
Let us give Brion-Szenes-Vergne (\cite{brionvergne},\cite{Szenesvergne}) formula  (slightly modified) for $\dim (\Sym_\lambda(N))$ in the form of an integral formula.
Let $g\in G$, and let $N^g$ be the fixed subspace of $N$ by the action of $g$.
If  the action of $G$ on $N^g$   is with generic finite stabilizer, we say that $g$ is a vertex of $\Phi$.
This is equivalent to say that the list $\Phi^g$ of  elements $\phi$ of $\Phi$ such that $g^{\phi}=1$ still spans $V$.
The set ${\CV}(\Phi)$ of vertices  is a finite subset of $G$.
The set of vertices ${\CV}(\Phi)$ is reduced to the identity element if and only if  the system $\Phi$ is unimodular.

We define $$A^{\C}(g)(X)=\det^{\C}_{N^g}\left(\frac{\exp(X)-1}{X}\right) \det_{N/N^g}^{\C}(1- g\exp(X))$$
and $$A^{\R}(g)(X)=\det^{\R}_{N^g}\left(\frac{\exp(X)-1}{X}\right) \det_{N/N^g}^{\R}(1- g\exp(X)).$$
Here $\det^{\C}$ is the complex determinant, while $\det^{\R}$ is the complex determinant.

Consider
$$Z(q,g)(X)=\frac{A^{\C}(g^{-1})(-X)}{A^{\R}(g)(qX)}$$
so
$$Z(q,g)(X)=\frac{\det^{\C}_{N^g}\left(\frac{1-\exp(-X)}{X}\right) \det_{N/N^g}^{\C}(1- g^{-1}\exp(-X))}
{\det^{\R}_{N^g}\left(\frac{\exp(q X)-1}{q X}\right) \det_{N/N^g}^{\R}(1- g\exp(q X))}.$$

For $q=1$,
$$Z(1,g)(X)=\frac{1}
{\det^{\C}_{N^g}\left(\frac{\exp(X)-1}{ X}\right) \det_{N/N^g}^{\C}(1- g\exp( X))}.$$

Remark that $Z(q,g)(X)$ is holomorphic at $q=0$.
Indeed $$Z(0,g)(X)=\frac{A^{\C}(g^{-1})(-X)}{ \det_{N/N^g}^{\R}(1- g)}.$$
Thus  $Z([q],g)(X)$ is  a series of analytic  functions of $X\in \g$.
Define the series $m([q],g)(y)$ of distributions on $\g^*$ by the formula
$$\la m([q],g),f\ra$$
$$ =\frac{1}{(2i\pi)^{\dim_\C N^g}}\int_{N^g} \left(\int_\g e^{\la X v,v\ra -\la dv,dv\ra }
Z([q],g)(X){\hat f} (X) dX\right).$$

Our set of admissible walls  is defined  to be $\CH:=\{p+H\}$ where $p$ is in the set
$F=\{-\phi_I\}$ and $H$ is a hyperplane generated by elements of $\Phi$.
Recall that the set $V_{reg}$ is the complement of the union of admissible walls.

\begin{theorem}\label{modifiedBZV}

 Write $m([q],g)(y)=\sum_{k=0}^{\infty} q^k m_k(g)(y)$.

$\bullet$  The distributions
 $m_k(g)(y)$  belong to the space $\CS$ of derivatives of piecewise polynomial generalized functions.

$\bullet$ If $g$ is not a vertex of $\Phi$,
the restriction of
$m_k(g)$ to $V_{reg}$ is equal to $0$.

$\bullet$ The restriction of
 $m_k(g)(y)$ to  each  connected component $\c$ of  $V_{reg}$ is a
polynomial function of $y$. It vanishes if $k\geq 2|\Phi|$.

$\bullet$ For any $\lambda\in \Lambda$, and any generic vector $\epsilon$,
$$\dim (\Sym_\lambda(N))=\sum_{g\in {\CV}(\Phi)} g^{-\lambda}\lim_{\epsilon} m([1],g)(\lambda).$$

\end{theorem}

The relation with the usual formulation of the value of the  Kostant partition function  $\dim (\Sym_\lambda(N))$
at  $\lambda$ via a sum of values  of splines functions obtained by differentiating $T(\Phi)$
is as follows.
The Fourier transform of the integral $I(X)=\int_{N^g}e^{\la X v,v\ra -\la dv,dv\ra }$ is the spline
$T(\Phi^g)$. Then, the Fourier transform of $A^\C(g^{-1})(-X)I(X)$ is the convolution of  $T(\Phi^g)$ with the box spline
$Box(-\Phi^g)$, followed by a series of translations. We obtain again a piecewise polynomial function on $V$.
Then to obtain the Fourier transform  of $\frac{1}{A^\R(g)([q]X)}A^\C(g^{-1})(-X)I(X)$, we apply an infinite series of constant coefficient differential operators to this piecewise polynomial function.

Remark that, in the last assertion  of Theorem  \ref{modifiedBZV}, $\epsilon $ is any generic vector. So
we tend to $\lambda$ coming from any direction.
So this theorem is not exactly Brion-Szenes-Vergne theorem (there it was required that $y$ tends to $\lambda$ along directions in the cone $\Cone(\Phi)$).
Here we have (as in the spirit of this article and of \cite{dpv1}) considered the representation of $G$ in $T^*N$ and written:
$$\prod_{\phi}\frac{1}{1-e^{i\phi}}=\prod_{\phi}\frac{(1-e^{-i\phi})}{(1-e^{-i\phi})(1-e^{i\phi})}$$
and expanded $$\Theta=\frac{1}{\prod_{\phi}(1-e^{-i\phi})(1-e^{i\phi})}$$
as a Fourier series supported in $\Cone(\Phi)$.
The theorem results then from the inversion formula for box splines \cite{DM1} (see also Theorem 2.29 in \cite{dpv2},
and \cite{vergnepoisson}).
Indeed, as the zonotope generated by $\Phi$ and $-\Phi$ contains $0$ in its interior, the values of the Fourier coefficients of $\Theta$
 are given by a piecewise analytic function on $V_{reg}$ which is continuous at any point of $\Lambda$

Furthermore if the set of vertices ${\CV}(\Phi)$ is reduced to $\{1\}$, then
the function $m([1],1)(y)$ gives by restriction to $V_{reg}$
a piecewise polynomial function which extends continuously to $V$ (see \cite{dpv2}, Remark 3.15).

Let $u\in U(N)$ be a unitary transformation of $N$ commuting with the action of $G$,
and let $R\in \End(N)$ be a complex endomorphism of $N$ commuting with the action of $G$ and with $u$.
The transformation $ue^R$ leaves stable the finite dimensional space  $\Sym_\lambda(N)$.
We consider the function $\lambda \mapsto \Tr_{\Sym_\lambda(N)} (u e^R)$.
Then, for $R$ small,
we can compute  the function $\lambda \mapsto \Tr_{\Sym_\lambda(N)} (u e^R)$
as a limit of a piecewise analytic functions on $V_{reg}$ (depending of $u,R$).
Let us give the formula.

%For each $g\in G$, we consider the space $N^{gu}$ fixed by the transformation $gu$. It is stable by the action of $G$ and $R$.
We define the set ${\CV}(\Phi,u)\subset G$ as the set of elements $g\in G$ such that the action of $g$ on $N^{gu}$ has a generic finite stabilizer.

Let
$$A^{\C}(g,u,R)(X)=\det^{\C}_{N^{gu}}\left(\frac{\exp(X+R)-1}{X+R}\right) \det_{N/N^{gu}}^{\C}(1- gu\exp(X+R))$$
and

$$A^{\R}(g,u,R)(X)=
\det^{\R}_{N^{gu}}\left(\frac{\exp( X+R)-1}{ X+R}\right) \det_{N/N^{gu}}^{\R}(1- gu\exp(X+R)).$$
Let
$$Z(q,g,u,R)(X)=\frac{A^{\C}(g^{-1},u^{-1},-R)(-X)}{A^{\R}(g,u,qR)(qX)}.$$
So

$$Z(q,g,u,R)(X)=\frac{\det^{\C}_{N^{gu}}\left(\frac{1-\exp(-(X+R))}{X+R}\right) \det_{N/N^{gu}}^{\C}(1- (gu)^{-1}\exp(-(X+R)))}
{\det^{\R}_{N^{gu}}\left(\frac{\exp(q (X+R))-1}{q (X+R)}\right) \det_{N/N^{gu}}^{\R}(1- gu\exp(q (X+R)))}.$$
We have
$$Z(1,g,u,R)(X)=\frac{1}
{\det^{\C}_{N^{gu}}\left(\frac{\exp( X+R)-1}{X+R}\right) \det_{N/N^{gu}}^{\C}(1- gu\exp(X+R))}.$$

Define the series $m([q],g,u,R)(y)$  of distributions on $\g^*$ by
the formula
$$\la m([q],g,u,R),f\ra$$
$$ =\frac{1}{(2i\pi)^{\dim_\C N^{gu}}}\int_{N^{gu}} \left(\int_\g e^{\la (X+R) v,v\ra -\la dv,dv\ra }
Z([q],g,u,R)(X){\hat f} (X) dX\right).$$

\begin{theorem}\label{maintheo}

Write $m([q],g,u,R)(y)=\sum_{k=0}^{\infty} q^k m_k(g,u,R)(y)$.

$\bullet$ If $g$ is not in  ${\CV}(\Phi,u)$,
the restriction of
$m_k(g,u,R)$ to $V_{reg}$ is equal to $0$.

$\bullet$  The restriction of
 $m_k(g,u,R)(y)$ to  each connected component $\c$ of  $V_{reg}$ is given by the restriction to $\c$ of  an analytic  function of $y\in V$.

 If $y\in \c$,
the series $$m([1],g,u,R)(y)=\sum_{k=0}^{\infty}  m_k(g,u,R)(y)$$  is convergent if $R$ is sufficiently small.

$\bullet$
If $R$ is nilpotent and $R^K=0$,
the restriction of
 $m_k(g,u,R)(y)$ to  each connected component $\c$ of  $V_{reg}$ is a
polynomial function of $y,R$, and vanishes if $k\geq K+2|\Phi|$.

$\bullet$ If $R$ is sufficiently small, then for any $\lambda\in \Lambda$, and any generic vector $\epsilon$,
$$\Tr_ {\Sym_\lambda(N)}(u e^R)=\sum_{g\in {\CV}(\Phi,u)} g^{-\lambda}
\lim_{\epsilon} m([1],g,u,R)(\lambda)$$
as a convergent series.

$\bullet$
If $R$ is nilpotent, then for any $\lambda\in \Lambda$, and any generic vector $\epsilon$,
$$\Tr_{\Sym_\lambda(N)}(u e^R)=\sum_{g\in {\CV}(\Phi,u)} g^{-\lambda}
\lim_{\epsilon} m([1],g,u,R)(\lambda)$$
as a finite sum.

\end{theorem}

Remark that when $X,R$ are small, we have the identities of meromorphic functions
$$Z(1,g,u,R)(X)\left(\int_{N^{gu}} e^{\la (X+R)v,v\ra } e^{-\la dv,dv\ra }\right)$$
$$=\det_{N^{gu}}^{\C}(\frac{X+R}{1-\exp(X+R)}) \det_{N/N^{gu}}^{\C}(\frac{1}{1-gu \exp(X+R)})(\frac{1}{\det^{\C}_{N^{gu}}(X+R)})$$
$$=\det_{N}^{\C}(\frac{1}{1-gu\exp(X+R)})=\Tr_{\Sym(N)}(gu e^{X+R}).$$

As $g\in G$, and $X\in \g$, this is
$$\sum_\lambda g^{\lambda} e^{i\la \lambda,X\ra } Tr_{\Sym_\lambda(N)}(ue^R)$$
which gives a  vague plausibility to the formula.

This theorem can be proved as in \cite{Szenesvergne}  by carefully reducing the computation  to a one dimensional computation.

\bigskip

We give some simple examples of the corresponding formula.

$\bullet$
Let $N$ be the complex $1$ dimensional space with action of $S^1$ by
$e^{i\theta}$.  Let
$u$ be a complex number of modulus $1$, and $r\in \C$.
Then we have, for $\lambda$  a non negative integer,
$$\Tr_{\Sym_\lambda(N)}(ue^{r})=u^{\lambda}e^{\lambda r}.$$

Let us compute our formula. Here we have just one vertex
$g=1/u$.
We compute the distribution  $m([1],g,u,r)(y)$,  for $g=1/u$, restricted to $V_{reg}$.
This is a piecewise analytic distribution on $\R$ given by the following
formula:

$$
m([1],g,u,r)(y)=
\begin{cases}
   e^{r y}\hspace{21mm}{\rm if}\ y\geq 0,\\
   (y+1)e^{r y}\hspace{9mm} {\rm if}\ -1\leq y\leq 0,\\
0\hspace{25mm}{\rm if}\ y\leq -1 .
\end{cases}
$$

So $m([1],g,u,r)(y)$ extends to a continuous function of $y$.
Furthermore we see that the limit  when $y\mapsto \lambda$ from left or right of $g^{\lambda}m([1],g,r)(y)$
is zero on all strictly negative integers, and is equal to $u^{\lambda}e^{\lambda y}$, when $\lambda$ is a non negative integer.

\bigskip

$\bullet$
Let $N$ be the complex two dimensional space with action of $S^1$ by
$(e^{i\theta},e^{2i\theta})$. Let $u\in S^1$ acting by  homothety  and let $R=0$.
Then $$Tr_{\Sym_\lambda(N)}(u)=
\begin{cases}\frac{u^{\lambda}}{1-u^{-1}}+\frac{u^{\lambda/2}}{1-u}\hspace{15mm}
{\rm if} \,\,\lambda {\rm \, is\, a \,non\, negative\, even\, integer},\\
\frac{u^{\lambda}}{1-u^{-1}}+\frac{u^{(\lambda+1)/2}}{1-u}\hspace{10mm}
{\rm if} \,\,\lambda {\rm\, is\, a\, positive\, odd\, integer},\\
0\hspace{35mm}{\rm if} \,\,\lambda {\rm \,is\, a \,strictly\, negative\, integer}.
\end{cases}
$$

Let us now compute our formula.
There are $3$ vertices, $g=1/u$ and $g_1,g_2$ the two elements such that $g_i^2=1/u$.

We compute the  distribution $m([1],1/u,u,0)(y)$, (restricted to $V_{reg}$).
This is a piecewise polynomial distribution on $\R$ given by the following
formula
$$m([1],1/u,u,0)(y)=
\begin{cases}
\frac{1}{(1-u^{-1})} \hspace{4cm} {\rm if}\,\, y\geq 0,\\
\frac{y}{(1-u^{-1})(1-u)}+\frac{u^2+3}{(1-u^{-1})(1-u)^2} \hspace{6mm} {\rm if}\,\, -1\leq y\leq 0,\\
\frac{1}{(1-u^{-1})^2} \hspace{39mm} {\rm if}\,\,  -2\leq y\leq -1,\\
-\frac{y}{(1-u^{-1})^2}-\frac{u-5}{(1-u)(1-u^{-1})^2}\hspace{9mm} {\rm if}\,\,   -3\leq y\leq -2,\\
0\hspace{51mm}  {\rm if}\,\,   y\leq -3.
\end{cases}$$

 Consider the vertices $g_1,g_2$.
We compute the distribution  $m([1],g_i,u,0)(y)$, (restricted to $V_{reg}$).
This is a piecewise polynomial distribution on $\R$ given by the following
formula:

$$m([1],g_i,u,0)(y)=\begin{cases}
\frac{1}{2(1-g_iu)}\hspace{50mm}{\rm if}\,\, y\geq 0,\\
\frac{1}{4}\frac{y}{(1-g_i u)(1-(g_iu)^{-1})}+\frac{1}{4}\frac{3u-3 g_i u+2}{(1-g_iu)^3}\hspace{10mm}{\rm if}\,\,  -1\leq y\leq 0,\\
 \frac{1}{4}\frac{y}{(1-g_i u)}-\frac{1}{4}\frac{3 g_i u-2}{(1-g_iu)^2}\hspace{29mm}{\rm if}\,\,  -2\leq y \leq -1,\\
 \frac{1}{4}\frac{y}{(1-g_i u)^2}-\frac{1}{2}\frac{2 g_i u-1}{(1-g_iu)^3}\hspace{28mm} {\rm if}\,\, -3\leq y\leq -2,\\
 0\hspace{61mm} {\rm if}\,\, y\leq -3.
 \end{cases}
 $$
One verifies that the limit  when $y\mapsto \lambda$ from left or right of $$\sum_{g\in {\CV}(\Phi,u)} g^{-\lambda} m([1],g,u,0)(y)$$
is equal to zero on all strictly negative integers, and is equal to $1$ when $\lambda=0$.
If $\lambda$ is a strictly  positive integer, we  have:

$$\sum_{g\in {\CV}(\Phi,u)} g^{-\lambda}m([1],g,u,0)(\lambda)=\frac{u^{\lambda}}{1-u^{-1}}+g_1^{-\lambda}\frac{1}{2(1-g_1u)}+ g_2^{-\lambda}\frac{1}{2(1-g_2u)}$$
 and this is indeed the right formula for $\Tr_{\Sym_\lambda(N)}(u).$

\section{The multiplicity index formula for a torus action}
Let $G$ be a torus acting on our manifold $M$, and let $V=\g^*$.
We assume that  $M$ is connected and that the generic infinitesimal stabilizer of the action of $G$ on $M$ is equal to $\{0\}$.
It is immediate to reduce to this case.
Furthermore, we  assume that our manifold $M$ can be embedded $G$-equivariantly in a vector space with a linear action of $G$.
This is to insure  that the set of stabilizers  of points of $M$ is a finite set of subgroups of $G$.

Let $\CS^1=\{\R U_a\}$ be the set of infinitesimal stabilizers for the action of $G$ on $M$ which are  of dimension $1$.
This is  a finite set of lines $\R U_a$ in $\g$.

\begin{definition}
We consider the following set of hyperplanes in $\g^*$
 $$\CH_0(M)=\{U_a^{\perp}; \R U_a\in \CS^1\}.$$
\end{definition}
If the action of $G$ on $M$ is locally free,  that is  if all stabilizers of points of $M$ are finite subgroups,  this set of hyperplanes is empty.

Let $\sigma$ be a transversally elliptic symbol on $M$.
Let $g\in G$ and  let $\Ch(g,\sigma)(X)$ be the twisted Chern character of $\sigma$.
If $\alpha\in H_G^*(M^g)$ is an equivariant cohomology class with polynomial coefficients,
 we can define a generalized function $h(g,\sigma,\alpha)$ on $\g^*$ by the formula
$$\int_{T^*M^g}^{\omega}\int_\g
\Ch(g,\sigma)(X)\alpha(X)\hat f(X)dX=\int_{\g^*} h(g,\sigma,\alpha)(y) f(y)dy$$ for any test
function $f$  on $\g^*$.

We say that $g\in G$ is a vertex of the action of $G$ on $M$ if there exists a point  $m\in M^g$ with  finite stabilizer.
Then the set $\CV(M)$ of vertices is a finite subset of $G$.
If all stabilizers are connected, the set $\CV(M)$ is reduced to the identity element.

 We have constructed the formal series
 $J([q],M^g) D([q],g, M/M^g)$ of equivariant classes with polynomial coefficients on $M^g$.
 Its inverse is also a series of equivariant classes with polynomial coefficients.

\begin{definition}
 Define the series of generalized functions $m([q],g,\sigma)$ of generalized functions  on $\g^*$ such that
$$\int_{T^*M^g}^{\omega}\int_{\g}  (2i\pi)^{-\dim M^g} \frac{{\rm Ch}(g,\sigma)(X)}{J([q],M^g)(X) D([q],g,M/M^g)(X)}\hat f(X)dX$$
$$=\int_{\g^*}m([q],g,\sigma)(y) f(y)dy$$ for any test function $f(y)$  on $\g^*$. \end{definition}

We write $m([q],g,\sigma)(y)=\sum_{k=0}^{\infty} q^k m_k(g,\sigma)(y)$, a series of generalized functions on $\g^*$.

The following lemmas and  theorem  will be proved in the next  subsections.

\begin{lemma}\label{lemmareg}
There exists a finite subset $\CH(\sigma)$ of affine walls such that
all distributions
$m_k(g,\sigma)$ are derivatives of piecewise  (with respect to $\CH(\sigma)$)  polynomial distributions.
Furthermore an affine hyperplane of the collection $\CH(\sigma)$ is parallel to an hyperplane in the family
$\CH_0(M)$.

\end{lemma}

Denote by  $V_{reg}$  the complement of the union of  affine walls in  $\CH(\sigma)$.
Let $\c$ be a connected component of $V_{reg}$.
It follows from Lemma \ref{lemmareg}  that
for each $k\geq 0$ and $g\in \CV(M)$, $m_k(g,\sigma)$ is given by a polynomial formula on each tope $\c$.

\begin{lemma}
 If $k\geq \dim M$, the distribution $m_k(g,\sigma)$ vanishes on $V_{reg}.$
\end{lemma}

Thus we can define the piecewise polynomial function $m(g,\sigma)$ on $V_{reg}$ by writing
$m(g,\sigma)=\sum_{k=0}^{\infty}m_k(g,\sigma)$.

\begin{theorem}\label{theo:formulaindexabelian}
For any $\lambda\in \Lambda$, and any generic vector $\epsilon$, we have
$$\mult_G(\sigma)(\lambda)=\sum_{g\in \mathcal V(M)} g^{-\lambda} \lim_{\epsilon}
m(g,\sigma)(\lambda).$$ \end{theorem}

This theorem is particularly nice when stabilizers of the action of $G$ on $M$ are connected.
In this case, there is only one vertex $g=1$, and the theorem reads

$$\mult_G(\sigma)(\lambda)=\frac{1}{(2i\pi)^{\dim M}}\int_{T^*M}^{\omega}\int_{\g} \frac{\Ch(\sigma)(X)}{J(M)[X]} e^{-i\la \lambda,X\ra}dX.$$
Here the notation $[X]$ in an equivariant cohomology class means that we have replaced the cohomology class  by its formal
expansion in homogeneous classes with polynomial coefficients. The value at $\lambda$ is obtained by limits from values in $V_{reg}$
 and the limit can be taken from any direction.
Furthermore,  in this case, the piecewise polynomial function $m(1,\sigma)$
on $V_{reg}$ extends continuously on $\g^*$. This is the function $m_G(\sigma)$ that we considered in the introduction.

This formula is particularly suggestive when $\lambda=0$.
Indeed, if $D$ is a transversally elliptic operator with principal symbol $\sigma$,
 the dimension of $[\Index_G(\sigma)]^G$ is the virtual dimension of the space of $G$-invariant solutions of $D$.
 It is given by
 $\int_G \Index_G(D)(g)dg$.
 On the other hand, for $g$ near $1$,
 $\Index_G(\exp X)=\frac{1}{(2i\pi)^{\dim M}}\int_{T^*M}^{\omega} \frac{\Ch(\sigma)(X)}{J(M)(X)}$.
 Our "miraculous formula" reads
  $$\int_G \Index_G(D)(g)dg=\frac{1}{(2i\pi)^{\dim M}}\int_{T^*M}^{\omega}\int_{\g} \frac{\Ch(\sigma)(X)}{J(M)([X])}dX.$$

Note the dichotomy in the formula: it is very important that $\Ch(\sigma)(X)$ is computed as a class with $C^{\infty}$ coefficients, and {\bf not} in the completion of $H_G^*(N)$.
 In contrast,
$\frac{1}{J(M)[X]}$ is in the completion of the ring $H_G^*(N)$ of equivariant characteristic classes with  polynomial coefficients.

\subsection{Generators for $G$-transversally elliptic symbols}
 Following Atiyah-Singer (see the monograph \cite{At}), we describe  a set of generators for
$K_G^0(T^*_GM)$.

Let $\s$ be a subalgebra of $\g$ arising as a stabilizer $\g_m$. The set of such subalgebras is finite.
We denote by $M((\s))$ the strata of $M$ consisting of elements $m$ such that
$\g_m=\s$.
Thus $M((\s))$ is a locally closed submanifold of $M$, eventually with several connected components denoted by $C_\s$.
 If $S$ is the subtorus of $G$ with Lie algebra $\s$,
the action of $G/S$ on $C_\s$ is infinitesimally free.
We identify a tubular neighborhood $U_\s$ of $C_\s$ in $M$ to the normal bundle $\CN_\s$ of $C_\s$ in $M$.
 Then $S$ acts fiberwise in $\CN_\s$, and the set of fixed points  $(\CN_\s)^S$ for the action of $S$ is reduced to the zero section.
  Thus the real vector bundle $\CN_\s$ is of even rank.
  Furthermore, we can choose a $G$-invariant complex structure $J$ on  the fibers of the bundle $\CN_\s$. Using $J$, consider $\CN_\s$ as a complex vector bundle. As explained in \cite{pep-vergnespin1}, we can then define a map
$At_J: K_G^0(T^*_G C_\s)\to K_G^0(T^*_GU_\s)$ using the lozenge product  of a $G/S$ transversally elliptic symbol
on $C_\s$
with the Atiyah symbol $at(\CN_\s)$.

The following theorem is proved in \cite{At}. See also \cite{pargen} for a detailed study of the
$R(G)$-module $K_G(T^*_GM)$.

\begin{theorem} \label{theogennerators}(Atiyah-Singer)
The module $K_G^0(T^*_G M)$ is generated  over $R(G)$
by the submodules  $At_J(K_G^0(T^*_G C_\s))$,
where $\s$ runs over
the infinitesimal stabilizers of the action of $G$ on $M$, and $C_\s$ over the connected components of $M((\s))$.
\end{theorem}

\subsection{The proof for an infinitesimally free action}\label{sub:inffree}

 Let us first make a remark on Fourier inversion.
Assume that $\Theta$ is a generalized function on $G$ such that
the  support $\CV$ of $\Theta$ is finite.
Then, for each $v\in \CV$ and  $X\in \g$ small,
$\Theta(v \exp X)=\int_{\g^*} e^{i\la y,X\ra }m_v(y) dy$
where $m_v(y)$ is a polynomial function on $\g^*$.
This equality is in the sense of generalized functions.
Then it is easy to verify (using for example Poisson summation formula)
 that, as a generalized function on  $G$,
$$\Theta(g)=\sum_{\lambda\in \Lambda}\left(\sum_{v\in \CV} v^{-\lambda}m_v(\lambda)\right) g^{\lambda}.$$
We thus see that the Fourier coefficients of $\Theta$  on $\hat G=\Lambda$
 are given by   the quasi polynomial function $\lambda\mapsto \sum_{v\in \CV} v^{-\lambda}m_v(\lambda)$.

\bigskip

Let $M$ be a manifold with an  infinitesimally free action of $G$.
The set $\CV(M)$ of vertices of the action is thus the set of elements $g\in G$ such that $M^g\neq \emptyset$.

The space $T^*_GM$ is a vector bundle over $M$, and using a $G$-invariant metric, we can write a direct sum decomposition of the bundle
$T^*M$ as
$T^*M=T^*_GM\oplus M\times \g^*$.
Thus we obtain a projection of $T^*M$ on $T^*_GM$ by
$(m,\xi+\eta)\mapsto \xi$, and
the moment map $\mu$ is  $\mu(m,\xi+\eta)=\eta$. Here
$m\in M$, $\xi\in (T^*_GM)_m$ and $\eta\in \g^*$  (and $T^*_GM$ is the set of zeroes of the moment map).

Recall the description of the equivariant cohomology of a manifold $M$ with infinitesimally free action.
We say that a differential form $\nu$ on $M$  is basic
if $\nu$ is invariant and if $\iota(v_X)\nu=0$ for all $X\in \g$. If $G$ acts infinitesimally freely on $M$,
the complex $\CA^{basic}(M)$ of basic forms is stable under the de Rham differential and
the equivariant cohomology of $M$ is the cohomology of the complex of basic forms.
This is due to H. Cartan in the polynomial case (see \cite{cartan}).  See \cite{duf-kum-ver} for extending the proof to the $C^{\infty}$-case.
In particular we have $H_G^k(M)=0$ if $k>\dim M-\dim G$. Thus any series
$\sum_{k=0}^{\infty} q^k \alpha_k$ of equivariant classes with polynomial coefficients is finite.

If $\CF$ is a $G$-equivariant vector bundle on $M$,
we can always choose a $G$-invariant basic connection $\nabla$ on $\CF$. That is,
the horizontal space on $\CF$ determined by
$\nabla$ contains the vector fields $v_X$ with $X\in \g$.
Thus the moment map determined by the connection is identically $0$, the equivariant curvature is just $R^{\CF}$, and the equivariant Chern character of $\CF$ is the basic form $\Tr(e^{R^\CF})$.

Let $\sigma$ be a transversally elliptic symbol on $M$.
When the action is free, the symbol $\sigma$ is simply the pull back of an elliptic symbol on the manifold $M/G$, and the distribution
 $\Index_G(\sigma)(g)$ is supported at $g=1$.
More generally, in the case of an infinitesimally free action,
it is easy to see  \cite{At} that the support of the generalized function $\Index_G(\sigma)(g)$ is contained in the set of the vertices $\CV(M)$.
Thus, by the remark above, the Fourier coefficients of $\Index_G(\sigma)(g)$
 are given by   a quasi polynomial function on $\Lambda$. The fact that this quasi polynomial function coincide with our
 formula will follow from the BV formula \cite{BV1},\cite{BV2}.

Provided $X$ is small enough, BV general formula asserts that, for any $g\in G$,
$$\Index_G(\sigma)(g\exp X)$$
$$=\int_{T^*M^g}\frac{1}{(2i\pi)^{\dim M^g}}e^{id_\g\omega(X)} \frac{\ChBV(g,\sigma)(X)}{J(M^g)(X)D(g,M/M^g)(X)}.$$

In BV formula, $\ChBV(g,\sigma)$ is a representative of the twisted equivariant Chern character
chosen to be  slowly increasing on the fibers of $T^*M$ and rapidly decreasing in the transverse directions $\xi$.
 The integral is shown to be convergent in the generalized sense.
 As the action is infinitesimally free, we can choose the differential form  $\ChBV(g,\sigma)$
 on $T^*M$ to be the reciproc image  of a basic form  on  $T^*_GM$ (by the projection $T^*M\to T^*_GM$),
and rapidly decreasing on the fibers of $T^*_GM$.
So in particular $\ChBV(g,\sigma)$ is independent of $X$.
The equivariant class  $J(M^g) D(g, M/M^g)$ can be represented by an invertible  basic form on $M$, so is independent of $X$.
  Thus
 $$I_g(X)= \int_{T^*M^g}\frac{1}{(2i\pi)^{\dim M^g}}e^{id_\g\omega(X)} \frac{\ChBV(g,\sigma)}{J(M^g)D(g,M/M^g)}$$
is defined for any $X\in \g$, in the generalized sense.
As $d_\g\omega(X)=\la \mu(m),X\ra+\Omega$, it is easy to see
 that  $I_g(X)$ is the Fourier transform of a polynomial distribution $h(g,\sigma)(y)$:
 $$\int_\g I_g(X)\hat f(X) dX=\int_{\g^*} h(g,\sigma)(y)f(y)dy$$
 for any test function $f$ on $\g^*$.
 We have to show that $h(g,\sigma)(y)=m([1],g,\sigma)(y).$

  In this work,  we have defined our distribution $m([1],g,\sigma)(y)$ using the infinitesimal index.
We have thus chosen
$\Ch(g,\sigma)$  to be represented by a differential form supported near $\supp(\sigma)$.
We can choose this differential form to be the reciproc image by  projection of $T^*M$ on $T^*_GM$ of a basic form with compact support on $T^*_GM$.
As in Remark \ref{defichBV}, the two choices $\ChBV(g,\sigma)$ and $\Ch(g,\sigma)$ differ by a  de Rham differential of a basic form
on $T^*_GM$ which is rapidly decreasing on the fibers of the bundle $T^*_GM\to M$.
We have  chosen
a function  $One_Z$ equal to $1$ in the neighborhood of $T^*_GM$. We might choose $One_Z(m,\xi+\eta)$
(where $\xi\in (T^*_GM)_m$ and $\eta\in \g^*$)  to be a function of $\eta$ compactly  supported on $\g^*$ and equal to
$1$ in a neighborhood of $0$.
Our series $$\frac{1}{J([q],M^g)(X) D([q],g,M/M^g)(X)}$$ is a finite series, and its value at $q=1$ is the basic form
$1/(J(M^g) D(g, M/M^g))$.
Thus we see as in  Lemma \ref{exinf} (see \cite{dpv2}, proof of Theorem 4.21)
that
$$ \int_{T^*M^g}\int_\g e^{id_\g\omega(X)} \frac{\ChBV(g,\sigma)}{J(M^g)D(g,M/M^g)}\hat f(X) dX$$
$$=\lim_{s\mapsto \infty} \int_{T^*M^g}\int_\g e^{is d_\g\omega(X)} One_Z
 \frac{\Ch(g,\sigma)}{J(M^g)D(g,M/M^g)}\hat f(X) dX.$$

With the help of the remark above, our theorem holds thus for
infinitesimally free action.

\subsection{The case of the Atiyah symbol on a Hermitian vector space}
Consider   a Hermitian vector space $N$.
Let $G$ be a torus acting on $N$, and we assume that
$N=\oplus_{\phi\in \Phi} N_\phi$ where the list $\Phi$ span a salient cone.
We use the notations of  Example \ref{defiat} and of Subsection \ref{subsecPW}.
The set $\CV(N)$ is  the set of vertices ${\CV}(\Phi)$.
The set of hyperplanes $\CH_0(N)$ is  the set of hyperplanes generated by elements of $\Phi.$

Then, the Atiyah symbol  $\sigma=at(N)$ is $G$-transversally elliptic, and
  $$\Index_G(\sigma)(g)=\Tr_{\Sym(N)}(g)=\sum_\lambda \dim(\Sym_\lambda(N)) g^{\lambda}.$$
Let us prove that our multiplicity index formula is true for the symbol $\sigma=at(N)$.

Let $g\in \CV(N)$.
The equivariant form $J(N^g)(X)$ is just equal to
$\det_{N^g}\frac{e^X-1}{X}$ while  the form
$D(g,N/N^g)(X)$ is $\det_{N/N^g}(1-g e^X)$ (here determinants are real determinants).
 Consider the series of distributions  $m([q],g,\sigma)(y)$ on $\g^*$ such that
$$\int_{T^*N^g}^\omega \int_\g \frac{\Ch(g,\sigma)(X)}
{J([q],N^g)(X)D([q],g,N/N^g)(X)}{\hat f}(X)dX$$
$$=\int_{\g^*} m([q],g,\sigma)(y)f(y)dy.$$

Let $\alpha$ be the $1$-form $Im(\la v, dv\ra)$ on the Hermitian space $N$.
We identify $T^*N$ with $N\oplus N$.
 Now we use the new variables $u=\xi-Jx$, $v=\xi+Jx$ on $N\oplus N$, and compute as in  Appendix 1 of \cite{BV1}.
Using the relation  between the Atiyah symbol, the Bott symbol and the Thom form, we
obtain
$$\int_{T^*N^g}^\omega \int_\g \frac{\Ch(g,\sigma)(X)}
{J([q],N^g)(X)D([q],g,N/N^g)(X)}{\hat f}(X)dX$$
$$=\lim_{s\mapsto \infty}\int_{N^g} \int_\g e^{i s d_\g\alpha(X)}
\frac{\det^{\C}_{N^g}\frac{(1-e^{-X})}{X}
\det^{\C}_{N/N^g}(1-g^{-1} e^{-X})}
{J([q],N^g)(X)D([q],g,N/N^g)(X)}.$$

We then use Formula \ref{modifiedBZV},
and we obtain our result.
Here the finite set of admissible walls is the set $p+H$ where $p$ varies in $F=\{-\phi_I\}$ and $H$ varies in $\CH_0(N).$

\subsection{The case of a vector bundle over a manifold with infinitesimally free action}\label{subvectorbundle}
Let $P$ be a manifold with an infinitesimally free action of a torus $G_1$.
Let us consider a $G_1$ equivariant Hermitian complex vector bundle $\CN\to P$, with complex structure $J$.
Let $G_2$ be a torus acting trivially on $P$ and fiberwise on $\CN$,
preserving the Hermitian structure, and commuting with the $G_1$ action.
We consider $G=G_1\times G_2$. So $\CN$ is a $G$-equivariant hermitian vector bundle over $P$.
We denote by $N$ the typical fiber of $\CN$.
This is a Hermitian vector space  provided with an action of  the torus $G_2$.

We write $\CN=\oplus_{i=1}^k \CN_{\phi_i}$ where $\phi_i$ are characters of $G_2$, and $\CN_{\phi}$ is the subbundle of $\CN$
where $G_2$ acts via  $g\cdot n=\phi(g)n$.
Each bundle $\CN_{\phi_i}$ is $G_1$-equivariant.
We assume that the list $\Phi=[\phi_1,\phi_2,\ldots, \phi_k]$
is contained in a half space.
Thus the bundle $$\Sym(\CN)=\oplus_{\lambda\in \hat G_2} \Sym_\lambda(\CN)$$ is a sum
of finite dimensional $G_1$ equivariant vector bundles over $P$.

Let $\sigma$ be a $G_1$ transversally elliptic symbol
on $P$.
Thus, for each $\lambda\in \hat G_2$, the symbol  $\sigma\otimes \Sym_{\lambda}(\CN)$
is $G_1$-transversally elliptic.

Let $at(\CN)$ be the Atiyah symbol.
We consider the $G$-transversally elliptic operator $At_J(\sigma)=\sigma\lozenge at(\CN)$
 on the total space of $\CN$.
We need to prove that our formula for the multiplicity index  holds for the transversally elliptic symbol
$At_J(\sigma)$.
 Recall (see \cite{pep-vergne:spinc}, Subsection 6.3)  that we have
$$\Index_G(At_J(\sigma))(g_1,g_2)=\sum_{\lambda_2\in \hat G_2}  \Index_{G_1}(\sigma\otimes \Sym_{\lambda_2}(\CN))(g_1) g_2^{\lambda_2}.$$
Thus we obtain, for $\lambda_1\in \hat G_1$ and $\lambda_2\in \hat G_2$,
$$\mult_G(At_J(\sigma))(\lambda_1,\lambda_2)=\mult_{G_1}(\sigma\otimes \Sym_{\lambda_2}(\CN))(\lambda_1).$$

Let us use the multiplicity formula obtained in Subsection \ref{sub:inffree} for the infinitesimally free action of $G_1$ on $P$.
If $v_1$ is a vertex of the action of $G$ on $P$, then $v_1$ produces a fiberwise transformation on $\CN$ restricted to  $P^{v_1}$.
We obtain
$$\mult_{G_1}(\sigma\otimes \Sym_{\lambda_2}(\CN))(\lambda_1)$$
$$=\sum_{v_1} v_1^{-\lambda_1} \int^{\omega}_{T^*P^{v_1}}\int_{ \g_1}
\frac{\Ch(v_1,\sigma\otimes \Sym_{\lambda_2}(\CN))}{J(P^{v_1})D(v_1,P/P^{v_1})} e^{-i\la \lambda_1,X_1\ra} dX_1.$$
We have written $J(P^{v_1})D(v_1,P/P^{v_1})$ for the invertible basic class representing the equivariant cohomology class
$J(P^{v_1})(X)D(v_1,P/P^{v_1})(X)$.
We have
$$\Ch(v_1,\sigma\otimes \Sym_{\lambda_2}(\CN))=\Ch(v_1,\sigma) \Ch(v_1,\Sym_{\lambda_2}(\CN)).$$
The value at $\lambda_1$ has to be computed by limit.

Let us be more explicit on $\Ch(v_1,\Sym_{\lambda_2}(\CN)).$
As the action of $G_1$ on  $P$ is infinitesimally free,  we can choose
a $G$-invariant hermitian connection $\nabla$ on $\CN$,  basic with respect to $G_1$.
Thus its equivariant curvature $R(X_1,X_2)$ is
$X_2+R$ where we still denote by $X_2$ the fiberwise action of  $X_2$ on $\CN$,
while $R$ is the curvature of $\nabla$.
Thus $R$  is an endomorphism of $\CN_p$ with  coefficients two-forms on $P^{v_1}$ and commuting with
the action of $v_1$ and $G_2$.  We thus have
$$\Ch(v_1, \Sym_{\lambda_2}(\CN))= \Tr_{\Sym_{\lambda_2}(\CN_p)}(v_1 e^{R}).$$

Consider now our formula. We use the notations of Subsection \ref{subsecPW}.
It is not difficult to see that the vertices $(v_1,v_2)$ for the action of $G=G_1\times G_2$ on $\CN$ are
the couples $v=(v_1,v_2)$ where $v_1\in {\CV}(P)$ while $v_2\in {\CV}(v_1,\Phi)$.
Write, in the sense of generalized functions,
$$\int^{\omega}_{T^*\CN^{v}}\int_{\g}
\frac{\Ch(v,At_J(\sigma))(X_1,X_2)}
{J([q],\CN^{v})(X_1,X_2)D([q],v,\CN/\CN^{v})(X_1,X_2)}
e^{-i\la y_1,X_1\ra } e^{-i\la y_2,X_2\ra}
 dX $$
$$=\sum_{k=0}^{\infty} m([q],v)(y_1,y_2).$$
Here $X=(X_1,X_2)\in \g=\g_1\times \g_2$ and $dX=dX_1dX_2$.
Our wanted formula is  $$\mult_G(At_J(\sigma))(\lambda_1,\lambda_2)=\sum_{(v_1,v_2)} v_1^{-\lambda_1}v_2^{-\lambda_2}
\lim_{\epsilon}m([1],v)(\lambda_1,\lambda_2).$$

We have, for $v=(v_1,v_2)$ acting by $u=v_1v_2$ on $\CN_p=N$,
$$\Ch(v,At_J(\sigma))(X_1,X_2)=\Ch(v_1,\sigma) \Ch(u,at(N))(X_2+R).$$
Similarly
$$J([q],\CN^{v})(X_1,X_2)=
J([q],P^{v_1})\det_{N^u}\left(\frac{\exp([q](X_2+R))-1}{[q](X_2+R)}\right),$$
$$D([q],v,\CN/\CN^{v})(X_1,X_2)=D([q],v_1,P/P^{v_1})\det_{N/N^{u}}(1-u \exp([q](X_2+R))).$$

Fix $v_1$ a vertex for the action of $G_1$ on $P$, and consider the integration over the fiber
$T^*N^{u}$  of the bundle $T^*\CN^{v}\to T^*P^{v_1}$.
Thus it is sufficient to prove that
$$\Tr_{\Sym_{\lambda_2}(\CN_p)}(v_1 e^{R})=$$
$$\sum_{v_2} v_2^{-\lambda_2}\int^{\omega}_{T^*N^u}
\int_{\g_2}\frac{\Ch(u,at(N))(X_2+R)}{J_{N^u}([q](X_2+R))D_{N/N^u}(u)([q](X_2+R)}e^{-i\la y_2,X_2\ra} dX_2.$$
Using the relation  between the Atiyah symbol, the Bott symbol and the Thom form,
we are reduced to Theorem \ref{maintheo}.
We remark that as $R$ is a matrix with valued $2$ forms on $P^{v_1}$, certainly the nilpotency index of
$R$ is smaller that the dimension of $P$, and our distributions vanish for $k\geq \dim \CN$.

Here our finite set of directions of hyperplanes is the set of the admissible hyperplanes of the form
$\g_1^*\oplus H$ where $H\subset \g_2^*$ is an admissible hyperplane for $\Phi$.
We see that that these hyperplaces are contained in $\CH_0(\CN)$.

We thus obtain our theorem for $At_J(\sigma)$.

\subsection{The general case for a torus}

Let $G$ be a torus with weight lattice $\Lambda$.
Let $M$ be a $G$-manifold and
let $\nu$ be a transversally elliptic symbol on $M$.
Let $a\in \Lambda$, and let $L_a$ be the corresponding $1$-dimensional representation of $G$.
If $\nu_a=\nu\otimes L_a$ is the twisted symbol of
$\nu$, we have
$\Index(\nu_a)(g)=g^a \Index(\nu)(g)$.
Remark that our multiplicity index formula is just
$\mult_G(\nu_a)(\lambda)=\mult_G(\nu)(\lambda-a)$.
Thus it is sufficient to prove our formula for generators  of $K_G^0(T^*_G M)$
over $R(G)$.

Let $\s$ be a subalgebra of $\g$ of the form $\g_m$, with $m\in M$.
Take a decomposition $G=G_1 \times G_2$, with $G_2$ the subtorus of $G$ with Lie algebra $\s$.
Then $G_2$ acts trivially on $M((\s))$, and $G_1$ infinitesimally freely.
From Subsection \ref{subvectorbundle}, our theorem holds for
transversally elliptic symbols $\nu=At_J(\sigma)$. Our proof of Theorem \ref{theo:formulaindexabelian}
follows from Theorem \ref{theogennerators}.

\subsection{An amusing example}\label{amusing}
\begin{figure}
\begin{center}
  \includegraphics[width= 3 in]{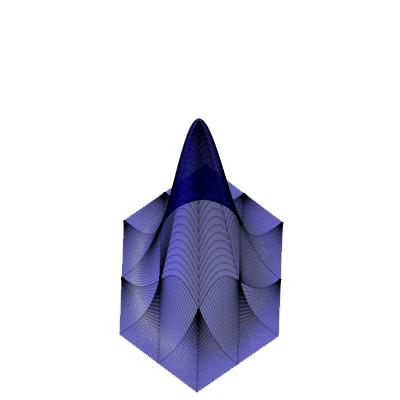}
\caption{The box spline for the root system $A_2$}
 \label{BA2}
\end{center}
\end{figure}

Let us consider the flag manifold  $M$ of $\C^3$ with the action of the adjoint torus $T$ in the adjoint group $K$ of $U(3)$.
 This is a case where all stabilizers are connected.
 Let $\Lambda$ be the root lattice of $A_2$. This is also identified to $\hat T$.
We consider the elliptic operator $\overline \partial$ on $M$, with index the trivial representation of $T$.
Let $\sigma$ its symbol. So the multiplicity index of $\sigma$  is the $\delta$ function on $\Lambda$. It takes value
$0$ at all points of $\Lambda$, except at $\lambda=0$ where the value is $1$.
We promised to obtain  (naturally) this multiplicity index as the restriction to $\Lambda$ of a continuous
spline on $\t^*$.
Let $\Phi$ be the root system $A_2$, consisting of $\pm\alpha, \pm \beta, \pm (\alpha+\beta)$.
Our formula is obtained by integration on $T^*(K/T).$
Using a similar computation as in the introduction, we see that
 $$B(X)=\int_{T^*(K/T)}\Ch(\sigma)(X)=\prod_{\phi\in \Phi} \frac{e^{i\la\phi,X\ra}-1}{i\la\phi,X\ra}.$$
This is the  Fourier transform of the box   spline $BA_2$ of the system $A_2$.
The box spline is supported on the convex hull of the points in $2\Phi$.
We draw in Figure \ref{BA2} the picture of the box spline $BA_2$.
It is  a locally polynomial measure, where each local polynomial is of degree $4$.
The value of $BA_2$ at $0$ is $\frac{1}{2}$.

\begin{figure}
\begin{center}
  \includegraphics[width=3 in]{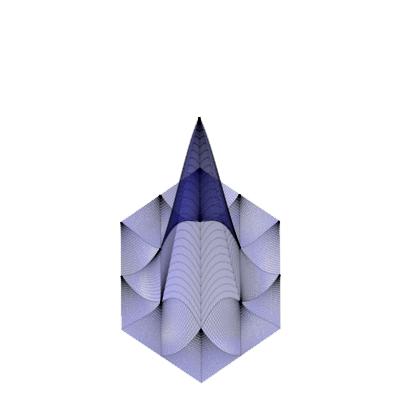}
\caption{The function $DBA_2$}
 \label{diffBA2}
\end{center}
\end{figure}

Now we consider the Taylor series $[B]$ up to order $4$ of the analytic function $B(X)$.
Formally, our function $m(\sigma)$ on $\t^*$  is the Fourier transform of
$B(X)/[B(X)]$ a function with Taylor series at $X=0$ equal to $1$ (however the Taylor series of $1/B(X)$ is not convergent!!).
 Thus $m(\sigma)$ is obtained by differentiating  $BA_2$  by the infinite series of differential operators $[B](i\partial)^{-1}$.
We draw the corresponding function $DBA_2$ on $\t^*_{reg}$ in Figure \ref{diffBA2}.

It is possible to verify directly that
the function $DBA_2$ is continuous and vanishes at all points of the root lattice $\Lambda$, except on $\lambda=0$, where its value is $1$.
This is fortunate, as this is the consequence of the inversion formula  for box splines (see for example \cite{vergnepoisson} for a proof).

\section{The case of a connected compact group}\label{secnonabelian}
We now give a formula for a general compact connected Lie group $G$.
We assume that the infinitesimal stabilizer of the action of the center of $G$  on $M$ is reduced to $0$.
We can always reduce easily to this case.

If $g \in G$, we denote by $G^g$ the centralizer of $g$ in $G$, and $\g^g$ its Lie algebra.
We say that $g$ is a vertex of the $G$-action if the infinitesimal stabilizer of the action of $G^g$ on $M^g$
is reduced to $0$.
Let $T$ be a maximal torus of $G$ with Lie algebra $\t$ and let $\CV(M)$ the set of vertices for the $T$-action on $M$.
Then $g$ is a vertex if and only if $g$ is conjugated to an element of $\CV(M)$.

We parameterize the set of irreducible representations of $G$ as follows.
We consider $\Lambda\subset \t^*$ the lattice of weights of $T$.
 We choose a system of positive roots $\Delta^+\subset \t^*$. Let
$\rho=\frac{1}{2}\sum_{\alpha>0}\alpha$, and $\t^*_{>0}$ be the positive (open) Weyl chamber. Let ${\tilde \Lambda}=\rho+\Lambda$
and denote by $\tilde \Lambda_{>0}={\tilde \Lambda}\cap \t^*_{>0}$.  For $\lambda\in \tilde \Lambda_{>0}$, we denote by
$V_\lambda$ the irreducible representation of $G$ of highest weight $\lambda-\rho$.

We consider the $G$-invariant function $d_G(y)$ on $\g^*$, where $d_G(y)$ is the volume of the coadjoint function $Gy$.
For $\lambda\in \tilde \Lambda_{>0}$, $d_G(\lambda)$ is the dimension of the representation $V_\lambda$.

Let $M$ be a $G$-manifold and $\sigma\in K_G^0(T^*_GM)$. We write
$$\Index_G(\sigma)=\sum_{\lambda\in \tilde \Lambda_{>0}} \mult_G(\sigma)(\lambda)  V_\lambda.$$

We now express  $\mult_G(\sigma)(\lambda)$
 in terms of the  equivariant Chern character of $\sigma.$

Let us consider the function $j_\g(X)=\det_\g \frac{e^{ad X}-1}{ad X}$.
Recall that $j_\g(X)$ has an analytic square root on $\g$.
We also defined $D_{\g/\g^g}(g)(X)= \det_{\g/\g^g}(1-g e^{ad X})$.
If $\rho$ is in the weight lattice $\Lambda$,
we choose the square root $D^{1/2}_{\g/\g^g}(g)(X)$ such that
$g^{-\rho}D^{1/2}_{\g/\g^g}(0)=i^{\dim \g^g/2}\prod_{\alpha>0}(1-g^{-\alpha})$.
Thus, for $\lambda\in \tilde \Lambda$,  $g^{-\lambda} D^{1/2}_{\g/\g^g}(g)(X)$ is well defined, without conditions on $\rho$.

Let $g\in G$.
 We have associated a series $J([q],M^g)(X) D([q],g,M/M^g)(X)$ of $G^g$ equivariant cohomology classes with polynomial coefficients on $M^g$.
We define the series of distributions
$m([q],g)$ on $(\g^g)^*$ by the formula
$$\la m([q],g),f\ra=$$
$$\int^{\omega}_{T^*M^g}\int_{\g^g }
\frac{j_{\g^g}^{1/2}(X)D^{1/2}_{\g/\g^g}(g)(X)\Ch(g,\sigma)(X)}{J([q],M^g)(X) D([q],g,M/M^g)(X)} {\hat f}(X) dX.$$

We write $m([q],g)(y)=\sum_{k=0}^{\infty}q^k m_k(g)(y)$.
Let us first give some indication on the nature of the generalized function $m_k(g)(y)$ on $(\g^g)^*$.
If $g\in T$, $\t$ is also a Cartan subalgebra of $\g^g$.
Thus, by $G^g$ invariance, we can restrict the function $d_{G^g}(y)m_k(g)(y)$ to  the interior $\t^*_{>0}$ of the Weyl chamber
(for $G$). We obtain a generalized function on $\t^*_{>0}$.
   In the course of the proof we will see that there exists a system $\CH$ of affine hyperplanes in $\t^*$ such that for any
   connected component $\c$ of the complement of $\CH$, the function
   $d_{G^g}(y)m_k(g)(y)$ is given by a polynomial function on $\c\cap \t^*_{>0}$.
   Furthermore, when $k$ is sufficiently large, $m_k(g)(y)$ restricts to $0$ on $\c\cap \t^*_>0$.
   Thus, for $\epsilon$ generic, and $\lambda\in \tilde \Lambda_{>0}$, $d_{G^g}(\lambda) \lim_{\epsilon} m([1],g)(\lambda)$
   is well defined.

\begin{theorem}
We have
$$\mult_G(\sigma)(\lambda)=\sum_{g\in \CV(M)} g^{-\lambda} d_{G^g}(\lambda) \lim_{\epsilon} m([1],g)(\lambda).$$
\end{theorem}

The theorem is particularly nice when stabilizers of the action of $G$ on $M$ are connected.
In this case, there is only one vertex $g=1$, and the theorem reads
$$\mult_G(\sigma)(\lambda)=\int^{\omega}_{T^*M}\int_{\g} j_\g(X)^{1/2}\frac{\Ch(\sigma)(X)}{J(M)[X]} e^{-i\la \lambda,X\ra}dX$$
a double integral formula reminiscent of Witten non abelian localization formula.

\begin{example}
Let $G$ acting on itself by left translations, and consider  $\sigma$ the $0$-symbol.
Its index is the $\delta$ distribution on $G$.
$$\Index_G(\sigma)=\sum_{\lambda\in \tilde \Lambda_{>0}} d_G(\lambda)  V_\lambda.$$

Let us compute our formula.
 The set of vertices is reduced to $1$, and the equivariant class $J(M)(X)$ is identically $1$.
 For any $s>0$, the distribution $I(X)=\int_{T^*G}e^{isd_\g\omega(X)}$ is easily seen to be the $\delta$ function on $\g$.

Thus our distribution is
$$\la m(\sigma),f\ra=\lim_{s\mapsto \infty} \int_{T^*G}\int_ \g
j_\g(X)^{1/2} e^{i sd_\g\omega(X)} {\hat f}(X) dX=\int_{\g^*} f(y) dy.$$
So $m([q],\sigma)(y)$ is identically equal to the function $1$ on $\g^*$, and we obtain our formula.
\end{example}

We now give the proof. It is based on  Theorem 5.21 in \cite{dpv2}
which relates the infinitesimal index  $\Infdex^{\omega}_G$ defined on $\CH_{G,c}^m(T_G^*M)$ and the infinitesimal index
$\Infdex_T^{\omega}$ defined on  $\CH_{T,c}^m(T^*_TM)$.

We may assume in the proof that $G$ is simply connected.
We recall that Atiyah-Singer associates to the $G$-transversally elliptic symbol $\sigma$ a $T$-transversally elliptic symbol
$A(\sigma)$ on $M$ such that
$$\mult_G(\sigma)(\lambda)=\mult_T(A(\sigma))(\lambda).$$

If $\sigma$ is itself $T$ transversally elliptic, we have, for $X\in \t$,
$$\Index_T(A(\sigma))(\exp X)=\prod_{\alpha>0} (e^{i\la \alpha/2,X\ra}-e^{-i\la\alpha/2,X\ra})\Index_G(\sigma)(\exp X).$$
More precisely, let $\mu: T^*M\to \g^*$ be the moment map. Write $\g=\t\oplus \q$.
Let $\mu^{\q}:T^*M\to \q^*$ be the projection of $\mu$ on $\q^*$.
We consider $T^*_TM$, relative to the action of $T$ on $M$. Thus $T^*_TM=T^*_GM\cap (\mu^\q)^{-1}(0)$.
The space $\q^*$ is provided with  the structure of a Hermitian vector space.
It thus have a Bott symbol $Bott(\q^*)\in K^0_T(\q^*)$ such that its restriction
to $\{0\}\in \q^*$ is $\prod_{\alpha>0}(1-t^{-\alpha})$.
Then,  as a $T$-transversally elliptic symbol,
$$A(\sigma)=\sigma \otimes (\mu^{{\mathfrak q}})^*Bott({\mathfrak q}^*)\otimes L_\rho.$$
Thus we compute $\mult_G(\sigma)(\lambda)$ by using our preceding theorem for $\mult_T(A(\sigma))(\lambda)$.

For $g\in \CV(M)$,
we need to compare the series  $m([q],g)(y)$ of distributions on  $\g^*$, and
 $n([q],g)(y)$ on  $\t^*$ associated respectively to
$\sigma$ and $A(\sigma)$

We do it for $g=1$.
By definition, $m([q],1)$ is the infinitesimal index of $\Theta(X)=j_\g^{1/2}(X)\frac{\Ch(\sigma)(X)}{J([q],M)(X)}$
and $n([q],1)$
is the infinitesimal index of $\theta(X)=\frac{\Ch(A(\sigma))(X)}{J([q],M)(X)}$.
Using the relation between the Bott symbol and the Thom class,
we see thus that $$\theta(X)=\frac{\Ch(\sigma)(X)}{J([q],M)(X)} (2i\pi)^{\dim \q}  (\mu^{{\mathfrak q}})^*(Thom(\q^*))(X).$$
By $G$-invariance, a $G$-invariant distribution restricts to $\t^*_{>0}$ (as a distribution).
\begin{lemma}\label{easydpv}
Let $P(X)\in \CH_{G,c}^m(T^*_GM)$ and let $$p(X)=(2i\pi)^{\dim \q} P(X) (\mu^{{\mathfrak q}})^*(Thom(\q^*))(X).$$
Then on $\t^*_{>0}$, we have
$$\Infdex^{\omega}_T(p)(y)=d_G(y)\Infdex_G^{\omega}(P)(y).$$
\end{lemma}

\begin{proof}
It is easy to see that this is true when $P(X)$ is in the cohomology with compact support of $T^*M$.
Then we can replace  $(2i\pi)^{\dim \q} Thom(\q^*)(X)$ by the function ${\bf Pf}(X)=\prod_{\alpha>0}\la \alpha,X\ra$.
The infinitesimal index is then just the Fourier transform on $\g$ of $I(X)=\int_{T^*M}P(X)$.
The formula follows from Harish-Chandra relation
between Fourier transform on $\g$ and $\t$ of $G$-invariant functions.

In the general case, we use Theorem 5.21 of \cite{dpv2}.
We have
$$\la \Infdex^{\omega}_G(P),f\ra =\la {\rm Ind}_{\t^*}^{\g^*}\nu,f\ra$$
with $\nu=\Infdex_T^{\omega}({\bf Pf}(X)p (X))$
and $$\la {\rm Ind}_{\t^*}^{\g^*}\nu,f\ra =\int_{\t^*} \nu(y)(\int_{{\mathfrak q^*}} f(y+q) dq)dy.$$
It is easy to see that this implies Lemma \ref{easydpv} by using Fourier transform.
 \end{proof}
The computation is similar for other $g\in \CV(M)$.


\begin{thebibliography}{99}










\bibitem{At} {Atiyah M. F.,} {Elliptic operators and compact groups.}  {\it Springer L.N.M.}, n. 401, 1974.\newline



\bibitem{Atiyah-Segal68} { Atiyah M.F., Segal G.B.},
{The index of elliptic operators II}, {\em Ann. Math.} {\bf 87} (1968),
531--545.\newline

\bibitem{Atiyah-Singer-1} {Atiyah M.F.,  Singer I.M.},
{The index of elliptic operators I},{\em Ann. Math.} {\bf 87} (1968),
484--530.\newline

%\bibitem{Atiyah-Singer-2} {\sc M.F. Atiyah}, {\sc I.M. Singer},
%{\em The index of elliptic operators III}, Ann. Math. {\bf 87},
%1968, p. 546-604.
%
%\bibitem{Atiyah-Singer-3} {\sc M.F. Atiyah}, {\sc I.M. Singer},
%{\em The index of elliptic operators IV}, Ann. Math. {\bf 93}, 1971,
%p. 139-141.





\bibitem{B-G-V} { Berline N., Getzler E.,   Vergne M.},
{\em Heat kernels and Dirac operators}, Grundlehren, vol. 298,
Springer, Berlin, 1991.\newline




\bibitem{BV82}{Berline N., Vergne M.,} {Classes caract\'eristiques \'equivariantes. Formule de localisation en cohomologie \'equivariante.}
{\it C.R. Acad. Sci.} {\bf 295} (1982), 539-541.\newline










\bibitem{BV2} {Berline N., Vergne M.,} { The Chern character of a transversally elliptic symbol and the equivariant index.}  {\it Invent. Math.} {\bf  124} (1996), 11--49. \newline

\bibitem{BV1} {Berline N., Vergne M.,}{ L'indice \'equivariant des op\'erateurs transversalement elliptiques.}
{\it Invent. Math.} {\bf 124} (1996), 51--101.\newline % % % % % % % %




\bibitem{brionvergne}
{Brion M., Vergne M.,} {Residues formulae, vector partition functions and lattice points in rational polytopes.}
 {\it Journal of the American Mathematical Society}
{\bf 10} (1997), 797--833.
\newline

\bibitem{DM1}{Dahmen W., Micchelli C., }
 {On the solution of certain systems of partial difference
  equations and linear dependence of translates of box splines.}
 {\it Trans. Amer. Math. Soc.}      {\bf 292}
  {(1985)}, 305--320.\newline


\bibitem{deboor} {De Boor C., H\"ollig K., Riemenschneider S.,}{Box Splines.}
{\em Applied Mathematical Sciences} {\bf 98}, Springer-Verlag New-York (1993).
\newline



\bibitem{dpv1}{De Concini C., Procesi C., Vergne M.,} {Vector partition functions  and index   of transversally elliptic operators.} {\it
    Transform. Groups} {\bf 15} {(2010)}, 775--811. \newline


\bibitem{dpv2}{De Concini C., Procesi C., Vergne M.,} {Box splines and the equivariant index theorem.}
    {\it  Journal of the Institute of Mathematics of Jussieu}
    {\bf 12}{(2013)}, 503--544.
arXiv :1012.1049 \newline


\bibitem{dpv3}{De Concini C., Procesi C., Vergne M.,} {The infinitesimal index.}
     {\it Journal of the Institute of Mathematics of Jussieu }{\bf 12}{(2013)}, 297--334.
      arXiv: 1003.3525
\newline



\bibitem{duf-kum-ver}{ Duflo M., Shrawan Kumar, Vergne M.,} {Sur la cohomologie \'equivariante des vari\'et\'es diff\'erentiables.}
{\em Ast\'erisque} {\bf 215} (1993).
\newline




\bibitem{Dui-Hec}{Duistermaat J.J, Heckman G.,} {On the variation in the cohomology of the symplectic form of the reduced phase space.}
{\it Invent. Math.}{\bf 69} (1982), 259--268.
\newline



\bibitem{Guiformal}{Guillemin V.},
{Reduced phase spaces and Riemann-Roch.}
Lie theory and geometry, 305–334,
Progr. Math., 123, Birkhäuser Boston, Boston, MA, 1994.
\newline



\bibitem{Guillemin-Sternberg82} { Guillemin V., Sternberg S.,} { Geometric quantization and multiplicities of group
representations.} {\em Invent. Math.} {\bf 67} (1982),  515--538.
\newline


\bibitem{cartan} { Guillemin V., Sternberg S.,} {Supersymmetry and equivariant De Rham theory. With an appendix containing two reprints by Henri Cartan.}{\em Mathematics Past and Present} Springer-Verlag, Berlin (1999).
\newline




\bibitem{Heckman82} {  Heckman G.J. ,} { Projections of orbits and asymptotic behavior of multiplicities
for compact connected Lie groups.} {\em Invent. Math.}{\bf 67} (1982),  333--356.
\newline


\bibitem{jk}{Jeffrey L. C., Kirwan F. C.,}
 {Localization and the quantization conjecture.} {\em Topology} {\bf 36} (1997), 647-693.
\newline




\bibitem{Meinrenken-Sjamaar} { Meinrenken E.,  Sjamaar R.},
{Singular reduction and quantization.} {\em Topology} {\bf 38} (1999),
 699--762.
\newline


\bibitem{pargen}{ Paradan P-E.,}
    {On the structure of $K_G(T_G M)$.}
     arXiv:1209.3852
\newline




\bibitem{par-ver4} { Paradan P-E, Vergne M.}, {The index of transversally elliptic operators.}
  {\it Ast\'erisque }{\bf 328} (2009),  297--338.
 \newline






\bibitem{pep-vergne:spinc}{ Paradan P-E, Vergne M.},
 {Multiplicities of equivariant Spinc Dirac operators.}
    arXiv:1411.7772 (2014).
\newline





\bibitem{pep-vergnespin1} { Paradan P-E, Vergne M.},
 {Witten non abelian localization for equivariant K-theory, and the [Q,R]=0 theorem.}
 arXiv:1504.07502  (2015).
\newline



\bibitem{sc} {Schoenberg I.J.,} {Cardinal Interpolation and Splines Functions II,}
{\em Journal of Approximation theory} {\bf 6}
(1972), 404-420.



\bibitem{Szenesvergne}{Szenes A., Vergne M.,} {Residue formulae for Vector partitions and Euler-MacLaurin sums,}
{\em Proceedings of FPSAC-01: Advances in Applied Mathematics}{\bf 30} (2003), 295--342.
\newline

%\bibitem{vergneicm}{Vergne M.},{Applications of Equivariant Cohomology.}
%{ICM congress Madrid} (2006)\newline


\bibitem{vergnepoisson}{Vergne M.,}{Poisson summation formula and Box splines,}
ArXiv:1302.6599.
\newline

\bibitem{Wittenequi} { Witten E.,} {Supersymmetry and Morse theory,}
{\em J.Diff. Geom.} {\bf 17} (1982), 661-692.
\newline


\bibitem{Wittennonabelian} { Witten E.}, Two dimensional gauge theories
revisited, {\em J. Geom. Phys.} {\bf 9}, 1992, p. 303-368.
\newline




\end{thebibliography}
\end{document}